\documentclass[aap,seceqn,MSNbibl,citesort,dvips]{arximspdf}

\doi{10.1214/09-AAP668}
\volume{20}
\issue{5}
\pubyear{2010}
\firstpage{1663}
\lastpage{1696}

\makeatletter

\newcommand{\iint}{\int\!\!\int}

\newtheorem{theorem}{Theorem}[section]
\newtheorem{corr}{Corollary}[section]
\newtheorem{lem}{Lemma}[section]

\def\bsuffix #1{#1}

\makeatother

\begin{document}
\begin{frontmatter}

\title{A time-dependent Poisson random field model for polymorphism
within and between two related biological species}
\runtitle{A time-dependent Poisson random field}

\begin{aug}
\author[A]{\fnms{Amei} \snm{Amei}\thanksref{t1}\ead[label=e1]{amei.amei@unlv.edu}} and
\author[B]{\fnms{Stanley} \snm{Sawyer}\corref{}\thanksref{t2}\ead[label=e2]{sawyer@math.wustl.edu}}
\runauthor{A. Amei and S. Sawyer}
\affiliation{University of Nevada, Las Vegas and Washington University
in St. Louis}
\address[A]{Department of Mathematical Sciences\\
University of Nevada, Las Vegas\\
Las Vegas, Nevada 89154\\
USA\\
\printead{e1}}
\address[B]{Department of Mathematics\\
Washington University in St. Louis\\
St. Louis, Missouri 63130\\
USA\\
\printead{e2}}
\end{aug}

\thankstext{t1}{Supported by Washington University
Dissertation Fellowship and NSF Grant DMS-01-07420.}

\thankstext{t2}{Supported in part by NSF Grant DMS-01-07420.}

\received{\smonth{3} \syear{2009}}
\revised{\smonth{10} \syear{2009}}

%
\begin{abstract}
We derive a Poisson random field model for population site polymorphisms
differences within and between two species that share a relatively recent
common ancestor. The model can be either equilibrium or time
inhomogeneous. We first consider a random field of Markov chains that
describes the fate of a set of individual mutations. This field is
approximated by a Poisson random field from which we can make inferences
about the amounts of mutation and selection that have occurred in the
history of observed aligned DNA sequences.
\end{abstract}

%

\end{frontmatter}

\section{Introduction}\label{sec1} A traditional goal of population genetics is to
understand the relationship between Darwinian selection and evolution.
A particular problem is to estimate the distribution of the fitness of
genes that become fixed in a natural population as a way of determining
whether evolution is going ``uphill'' or ``downhill'' in that population
(both are possible; \cite{rNatMod,rHartlPrimer,rHartlClark,rLewontin,rLiWH}).

One approach is to use the numbers of site polymorphisms between and
within a pair of closely related species at a genetic locus (McDonald and
Kreitman \cite{rMcDKr}; see also
\cite
{r14,rCaicedo,rDurr,rEyreWalker,rHartlPrimer,rHartlClark,rKeight1,rKeight2,r66}).
Sawyer and Hartl \cite{rSSDH} developed a Poisson random field (PRF) model for
these counts that can be used to estimate the amounts of selection that
are occurring at individual loci (\cite{r11,rDirSel,rDHMo,rLiWH,rSSGo};
see also \cite{rSewallWr}). Multilocus
extensions of the model can be used for large numbers of different loci
\cite{rBaines,rBoyko,rNatMod,rSSKul,rSSPar,rWilmson}. The models
assume a high level of recombination between nucleotides and are well
suited for the analysis of polymorphism and divergence at multiple loci
distributed across a genome \cite{rNatMod}. Other authors have extended
the basic PRF model to more general biological settings
\mbox{\cite{r24,r44,rSSPar,r71,r72,rWilmson}} and have used numerical
simulations to study
alternative models and the effects of deviations from model assumptions
\cite{rHaley,rBoyko,rDirSel,rZhuLan}. Williamson et al.
\cite{rWilmson} used a time-dependent PRF model to estimate the time since
a hypothetical sudden change in human population size. The purpose of this
paper is to provide a rigorous derivation of a time-inhomogeneous PRF
model by approximating discrete time Markov chains by diffusion processes
and, so that the model can be applied to data, to derive the corresponding
sampling formulas for aligned DNA sequences.

First, a brief introduction to genetics and population genetics. The
genetic material or DNA of most living creatures is located in one or more
chromosomes. A~chromosome can be thought of as a long string of letters
from the alphabet A, C, G, T, where each letter corresponds to a specific
nucleotide. Most higher plants and animals are \textit{diploid}, which means
that their DNA is arranged in pairs of chromosomes. Each of these paired
chromosomes typically has double-stranded DNA, amounting to four DNA
strands for each chromosome pair. The nucleotides on the two DNA strands
of each chromosome are mirrored as a way of providing redundant
information. Creatures with nonpaired chromosomes (such as bacteria and
viruses) are called \textit{haploid}. A \textit{gene} or \textit
{genetic locus}
is a
relatively-short segment of a chromosome that affects a particular trait
or set of traits. \textit{Random mating} for a population of $M$ diploid
individuals, along with Mendel's laws for the offspring of two diploid
individuals, can be statistically approximated by random sampling with
replacement of $N=2M$ haploid individuals to form the next generation.
That is, random mating of $M$ diploid individuals is approximated by
random mating of their $N=2M$ haploid genomes.

Proteins are built from peptides, which are strings of amino acids.
Peptide strings are created by reading a sequence of \textit{codons}, which
are consecutive triples of nucleotides, from one or more subsegments of a
gene. One amino acid is added to the peptide for each codon of DNA\null.
(The amino acids and codons have no chemical similarity to one another.)
There are around twenty different amino acids as opposed to $4^{3}=64$
possible codons on one DNA strand, so that the mapping of codons to amino
acids cannot be one to one. About half of amino acids are encoded by four
different codons that vary in the third codon position. Most of the
remaining amino acids are encoded by two different codons. A mutation at
one of the three sites in a codon position is called a \textit{silent}
mutation if the resulting new codon encodes the same amino acid. A site
mutation that results in a different amino acid is called a
\textit{replacement} mutation. Most mutations at first and second position
sites in codons are replacement, while mutations at the third position
can be either silent or replacement. The majority of replacement mutations
severely damage the gene product and may be lethal to the host, but a
substantial number have mild to moderate effects and can be advantageous.

A DNA site is \textit{polymorphic} in a population if there is more
than one
nucleotide at that site. If there is only one nucleotide at that site
in a
population, the site is \textit{monomorphic} and the population is
\textit{fixed}
at that nucleotide.

The \textit{fitness} of a gene is defined as the expected relative
number of
surviving offspring (genes) in the next generation that are descendants of
that gene, assuming constant lifetimes for individuals. Suppose that we
have $N$ haploid individuals that are of one of two types $A$ or $a$ at a
particular site or gene, where $a$ is the mutant type and $A$ is the
pre-existing or ``wild type.'' Assume that $a$ has the fitness coefficient
$\omega_N(a)=1+\sigma_N$ relative to the pre-existing $A$ with
$\omega_N(A)=1$.
If $\sigma_N$ is scaled as $\sigma_N\sim\gamma/N$ for a large haploid
population size $N$, then $\gamma$ is the \textit{selection
coefficient}
of the
mutant. The mutation is favorable if $\gamma>0$, detrimental if
$\gamma<0$ and
neutral if $\gamma=0$.

\section{Main results and discussion}\label{sec2} We first derive two results that lead
to Poisson random fields (see below) for the distribution of site
polymorphisms in a limiting infinitely large random-mating population.

The population results cannot be directly applied to data, both because
we cannot sample the entire population and also because the results imply
that there are an infinitely large number of site polymorphisms in each
gene, most of which have very small population frequencies (see below).
However, we can use the population results to derive the distribution of
sample statistics for aligned DNA sequences arising from two related but
not too distantly related species. The sample statistics will turn out to
be independent Poisson with means depending on population parameters for
mutation and selection rates. These, in turn, lead to likelihood methods
for the estimation of mutation, selection and divergence time parameters.

The proofs of the population results are deferred to Sections \ref{sec4}--\ref{sec8}.
Mathematically, the key steps in the proofs of both of the main population
results depend on Trotter's theorem for the dual process of a diffusion
process (Section \ref{sec6}). The sampling formulas are discussed in Sections \ref{sec22}
and \ref{sec23} with the details in Section~\ref{sec3}.\looseness=1

\subsection{Population formulas}\label{sec21}
The first result describes the
limiting distribution of the population proportions of mutant nucleotides
at polymorphic sites at a single genetic locus in a large population.
Under the assumptions of Section \ref{sec1}, this distribution is a Poisson random
field (PRF) \cite{rKingman} on $(0,1)$, where each point is a population
frequency ratio at a different site (see Theorem \ref{theo21} below). The second
result (Theorem \ref{theo22}) describes the expected number of polymorphic sites
that have become fixed in the population at the mutant nucleotide since
time 0. This random variable is also Poisson and is independent of the PRF
derived in Theorem \ref{theo21}.

Silent and replacement polymorphic sites are modeled separately and
provide different information (see below).

We begin with a model for the frequency of a mutant nucleotide at a
single site in a finite population of haploid size $N$. Assume that
mutation is sufficiently rare at the site level that repeat mutations at
the same site can be neglected. We use Moran's second model \cite{rMoran}
for which, at each step, an individual is randomly chosen from the
population with equal probabilities and replaced by a copy of a second
individual (perhaps the same as the first) chosen with probability
proportional to the fitness. Then the number $X_k$ of individuals in the
population that have the mutant nucleotide at a particular site is a
Markov chain on $S_N=\{0,1,2,\ldots,N\}$ with transition function
%
\begin{eqnarray}\label{equ21}
p_N(i,i+1) &=& \frac{(1+\sigma_N)i/N(1-i/N)}
{1+\sigma_N i/N},\nonumber\\
p_N(i,i-1) &=& \frac{i/N(1-i/N)}
{1+\sigma_N i/N},\\
p_N(i,i) &=& 1-p_N(i,i-1)-p_N(i,i+1) \nonumber
\end{eqnarray}
for $0<i<N$ with $p_N(0,0)=p_N(N,N)=1$, where $1+\sigma_N$ is the
fitness of
the mutant nucleotide $a$ with respect to $A$. A generation corresponds to
$N$ time steps. The states $0$ and $N$ are traps corresponding to the loss
of all mutant ($a$-type) or original ($A$-type) nucleotides, respectively.
We write $X^N_k$ instead of $X_k$ when we want to emphasize the size of
the state space $N$.

We assume that, at time 0, there are $M_0$ different sites that are
population polymorphic with both mutant and nonmutant nucleotides. There
are mutations at $M_r$ new sites at times $r=1,2,\ldots.$ The $M_0$
and $M_r$ are independent Poisson with $E(M_r)=\mu_N$ and typically
$E(M_r)\ll E(M_0)$ ($r\ge1$). All mutations are assumed to occur at new
sites. Given $M_0$ and $M_r$, the numbers of mutant nucleotides at these
sites at times $k\ge0$ in the future (for the initial polymorphic sites)
and at times $k\ge r$ (for the new mutations) are given by independent
Markov chains $X_{1,a,k}$ ($1\le a\le M_0$) and $X_{2,b,r,k}$ ($1\le
b\le
M_r$), each of which has the transition function (\ref{equ21}). By assumption,
$X_{2,b,r,r}=1$ for $r\ge1$ and $0\le X_{1,a,k}, X_{2,b,r,k}\le N$.

The number of sites at which there are $j$ mutant nucleotides at time $k$
($1\le j\le N, k\ge0$) is
%
\begin{equation}\label{equ22}
N_k(j) =\#\{a\dvtx X_{1,a,k}=j\}
+ \#\{(b,r)\dvtx X_{2,b,r,k}=j\},
\end{equation}
where $\#$ means cardinality, $1\le a\le M_0$, $1\le b\le M_r$
and $1\le r\le k$. Thus
%
\begin{equation}\label{equ23}
\sum_{j=1}^N f\biggl(\frac jN\biggr) N_k(j) =
\sum_{a=1}^{M_0} f\biggl(\frac{X_{1,a,k}}N\biggr)
+\sum_{r=1}^k \sum_{b=1}^{M_r}
f\biggl(\frac{X_{2,b,r,k}}N\biggr)
\end{equation}
for functions $f(x)$ on $[0,1]$. The expected value of $N_k(j)$ in (\ref{equ22})
is
%
\begin{equation}\label{equ24}
E(N_k(j))
=\sum_{i=1}^{N-1} \omega^N_i p_N^{k}(i,j)
+\mu_N \sum_{r=1}^{k} p_N^{k-r}(1,j),
\end{equation}
where $\omega^N_i = E(N_0(i))$ and $p_N^k(i,j)$ is the
$k$th matrix
power of $p_N(i,j)$. The random variables $N_k(j)$ are somewhat
unintuitive in that each variable depends on information from many
different polymorphic sites, rather than one, but turn out to have good
statistical properties.

We assume that the random variables $M_r$ and $N_0(i)$ are independent and
Poisson distributed. Then, by a slight extension of what is called
Bartlett's theorem~\cite{rKingman}, this implies that, for each fixed
$k>0$, the counts $N_k(i)$ ($i=1,2,\ldots$) are independent Poisson and
thus define a Poisson random field (PRF) on $I_N=\{1/N,2/N,\ldots,1\}$.
If the mean measures defined by (\ref{equ24}) for $y=j/N$ and $k=k_N$ converge
weakly as $N\to\infty$ on compact subsets of $(0,1]$, then $N_{k_N}(i)$
for $x=i/N$ converge weakly in the same sense to a PRF on $(0,1]$. The
points of the limiting PRF will be a countably-infinite set of diffusion
processes on $(0,1]$ each fixed at a particular time $t$.

Each of the Markov chains $X_{1,a,k},X_{2,b,r,k}$ can be approximated for
large $N$ by a diffusion process $X_t$ on $(0,1)$ with time scaled as
$t\sim k/N^2$ (Section \ref{sec4}). The diffusion process $X_t$ has traps at the
endpoints $0,1$, similarly to the Markov chains. Specifically, $X_t$ is
the diffusion process on $(0,1)$ generated by the differential operator
%
\begin{equation}\label{equ25}
L_x =x(1-x)\,\frac{d^2}{dx^2} +\gamma x(1-x)\,\frac{d}{dx},
\end{equation}
where $\gamma=\lim_{N\to\infty}N\sigma_N$ for $\sigma_N$
in (\ref{equ21}). We write $L_x$
in the Feller form
%
\begin{equation}\label{equ26}
L_x=\frac{d}{m(dx)}\,\frac{d}{s(dx)},
\end{equation}
where $m(dx)=m'(x)\,dx$ and $s(dx)=s'(x)\,dx$ for
%
\begin{equation}\label{equ27}
s(x)=\frac{1-e^{-\gamma x}}{\gamma} \quad\mbox{and}\quad
m(dx)=\frac{e^{\gamma x}}{x(1-x)}\,dx.
\end{equation}
The function $s(x)$ and measure $m(dx)$ are called the \textit{scale} and
\textit{speed measure} of $X_t$, respectively. The diffusion process
$X_t$ has
a smooth symmetric transition density $p(t,x,y)=p(t,y,x)$ with respect
to $m(dx)$ such that
%
\begin{eqnarray}\label{equ28}
E_x(f(X_t))
&=& Q_tf(x) =E\bigl(f(X_t)\mid X_0=x\bigr) \nonumber\\[-8pt]\\[-8pt]
&=& \int_0^1 p(t,x,y)f(y)m(dy) \nonumber
\end{eqnarray}
for $f\in C[0,1]$ with $f(0)=f(1)=0$.

Fix $t>0$ and assume that $\sigma_N,\mu_N,k_N$ are scaled so that
%
\begin{equation}\label{equ29}
N\sigma_N\to\gamma,\qquad N\mu_N\to\theta\quad
\mbox{and}\quad k_N/N^2\to t \qquad\mbox{as $N\to\infty$}.
\end{equation}
As mentioned earlier, each step in the Moran model (\ref{equ21}) puts one
individual at risk, so that each individual is put at risk on the average
every $N$ time steps. Thus $N\sigma_N\to\gamma$ and $k_N\sim tN^2$
mean that
selection at intensity $\gamma$ is applied to each individual in the
diffusion time scale $t$.

The relation $N\mu_N\to\theta$ implies that the rate of arrival of new
mutations in the diffusion time scale is $(N^2)\mu_N\sim N\theta\to
\infty$.
However, the new mutant Markov chains begin at $X_{2,b,r,r}=1$ which
corresponds to states $1/N\to0$ for the approximation diffusion. The
limiting processes $X_t$ have a trap at $x=0$, so that one might expect
that only $O(1/N)$ of the new mutant Markov chains survive the first few
generations. This would suggest that only of order
$O((1/N)(N^2\mu_N))=O(\theta)$ of the new mutants would
survive, so
that the rate $N\mu_N\to\theta$ in (\ref{equ29}) is not paradoxical.

The equilibrium distribution of the limiting PRF of diffusion processes
is
%
\begin{equation}\label{equ210}
\mu_{\theta,\gamma}(dx) =
\theta\frac{s(1)-s(x)}{s(1)-s(0)} m(dx)
=
\frac{s(1)-s(x)}{x(1-x)} \frac{\theta e^{\gamma x}}{s(1)} \,dx
\end{equation}
(\cite{rSSDH,rSewallWr}; see also Corollary \ref{corr21} below). In particular, the
equilibrium mean density has a $1/x$ singularity at $x=0$ in a large
population but is bounded at $x=1$. This means that, in the limit, there
are an infinite number of sites that have polymorphic mutant alleles with
population proportions $p_i>0$. By properties of Poisson random fields,
\[
E\biggl(\sum_i p_i\biggr) =\int_0^1 x \mu_{\theta,\gamma}(dx)
=\frac{\theta}{s(1)} \int_0^1 \frac{s(1)-s(x)}{1-x}e^{\gamma
x} \,dx
<\infty.
\]
Thus the vast majority of the $p_i$ are too small for the mutant
nucleotide to show up in finite samples. However, this does show that we
need to allow for initial limiting PRF mean measures that are not
normalizable at zero.

Specifically, for the initial distribution of polymorphic sites, we
assume that there exists a Borel measure $\nu(dx)$ on $(0,1)$ such that
$\int_0^1 x \nu(dx)<\infty$ and
%
\begin{equation}\label{equ211}
\lim_{N\to\infty} \sum_{j=1}^{N-1}
g \biggl(\frac jN \biggr) \frac jN \omega^N_j
=\int_0^1 g(y) y \nu(dy)
\end{equation}
for all $g\in C[0,1]$ for $\omega^N_j=E(N_0(j))$ in (\ref{equ24}). Since $s'(0)=1$
in (\ref{equ27}), this is equivalent to
%
\begin{equation}\label{equ212}
\lim_{N\to\infty} \sum_{j=1}^{N-1}
g \biggl(\frac jN \biggr) s \biggl(\frac jN \biggr) \omega^N_j
=\int_0^1 g(y)s(y)\nu(dy).
\end{equation}
Our first result describes the limiting PRF on $(0,1)$ which in turn
describes the population distribution of polymorphic mutant sites. The
proof is deferred to Section \ref{sec7}.
\begin{theorem}\label{theo21}
Assume that $\sigma_N,\theta_N,k_N$ satisfies (\ref{equ29})
and that
$N_k(i)$ defined in (\ref{equ22}) satisfies (\ref{equ211}) for $\omega^N_j=E(N_0(j))$. Then
for $Q_tf(x)$ in (\ref{equ28})
%
\begin{eqnarray}\label{equ213}\qquad
&&
\lim_{N\to\infty} E\Biggl(\sum_{i=1}^{N-1} f\biggl(\frac iN\biggr) N_{k_N}(i)\Biggr)
\nonumber\\[-8pt]\\[-8pt]
&&\qquad=
\int_0^1 Q_tf(x) \nu(dx) +
\theta\int_0^1 \frac{s(1)-s(x)}{s(1)-s(0)}
\bigl(f(x) - Q_tf(x)\bigr) m(dx) \nonumber
\end{eqnarray}
for any $f\in C[0,1]$ with $f(0)=f(1)=0$ such that $g(x)=f(x)/x$ for $x>0$
extends to a continuous function on $[0,1]$.
\end{theorem}

The limiting PRF mean density in (\ref{equ213}) is $g(t,\theta,\gamma
,y)m(dy)$ where
%
\begin{eqnarray}\label{equ214}
g(t,\theta,\gamma,y) &=& \int_0^1 p(t,x,y)\nu(dx) +
\theta\frac{s(1)-s(y)}{s(1)-s(0)}\nonumber\\[-8pt]\\[-8pt]
&&{}
-\theta\int_0^1
\frac{s(1)-s(x)}{s(1)-s(0)} p(t,x,y)m(dx). \nonumber
\end{eqnarray}
The first terms on the right in (\ref{equ213}) and (\ref{equ214}) are transient terms
that are due to the initial (or ``legacy'') polymorphisms at time $t=0$.
The remaining terms are due to new mutations that were introduced at times
$t>0$.

If $\nu(dx)$ in (\ref{equ211}) and (\ref{equ214}) is the measure $\mu_{\theta,\gamma}(dx)$
in (\ref{equ210}), then the right-hand side of (\ref{equ213}) is identically
$\int_0^1 f(x) \mu_{\theta,\gamma}(dx)$, so that $\mu_{\theta
,\gamma}(dx)$ is an
equilibrium measure. The following corollary shows that
$\mu_{\theta,\gamma}(dx)$ is an asymptotic measure as well (Section \ref{sec73}).
\begin{corr}\label{corr21} Let $G_t(f)$ be the right-hand side of (\ref{equ213}) for any function
$f(x)$ satisfying the conditions of Theorem \ref{theo21}. Then
%
\begin{equation}\label{equ215}\quad
\lim_{t\to\infty} G_t(f) =
\theta\int_0^1 \frac{s(1)-s(x)}{s(1)-s(0)} f(x) m(dx)
=\int_0^1 f(x) \mu_{\theta,\gamma}(dx).
\end{equation}
\end{corr}

Theorem \ref{theo21} implicitly assumes that all new replacement mutants at a given
genetic locus have the same selection coefficient $\gamma$. The model
can be
extended to allow within-locus random distributions of selective effects
of new replacement mutants, for example, Gaussian
\cite{rBaines,rSSKul,rSSPar}, exponential \cite{rHaley} or
gamma \cite
{rBoyko}.

The second main population result describes the limiting expected number
of mutant sites that have become fixed in the population at the mutant
nucleotide by time $t>0$ (i.e., such that the nonmutant or ``wild
type'' nucleotide has been lost). The limit can be expressed in terms of
the hitting times $T_a=\min\{{t\dvtx X_t=a}\}$ for the limiting diffusion
process $X_t$ and for its dual process $\widetilde{X}_t$, but can also
be expressed
in terms of the transition density of $X_t$. The proof is deferred to
Section~\ref{sec8}.
\begin{theorem}\label{theo22}
Under the conditions of Theorem \ref{theo21}, the asymptotic expected number of
mutant sites that have become fixed in the population at the mutant
nucleotide by time $t$ is
%
\begin{eqnarray}\label{equ216}
&&\lim_{N\to\infty} E (N_{k_N}(N)) \nonumber\\
&&\qquad =\lim_{N\to\infty}
\Biggl(\sum_{i=1}^{N-1} \omega^N_i p_N^{k_N}(i,N)
+\mu_N \sum_{r=1}^{k_N} p_N^{k_N-r}(1,N) \Biggr) \\
&&\qquad =
\int_0^1 P_x(T_1\le t) \nu(dx)
+\frac\theta{s(1)}
\int_0^t \widetilde{P}_0(T_1\le u) \,du, \nonumber
\end{eqnarray}
where $\widetilde{P}$ in (\ref{equ216}) refers to the diffusion process $X_t$
conditioned
on $T_1<T_0$, which has an entrance boundary at $x=0$ (Section \ref{sec6}).
\end{theorem}

The right-hand side of (\ref{equ216}) can also be written (Section \ref{sec82})
%
\begin{eqnarray}\label{equ217}
&&\frac1{s(1)} \biggl(\int_0^1 s(x)\nu(dx) -
\int_0^1\int_0^1 p(t,x,y)s(y)m(dy) \nu(dx) \nonumber\\[-8pt]\\[-8pt]
&&\hspace*{37.7pt}\qquad{} +\theta t -\theta\int_0^t\int_0^1
q(u,0+,y)s(y)^2m(dy) \,du\biggr), \nonumber
\end{eqnarray}
where $q(t,x,y)$ is the transition density of the dual-process diffusion
(Section \ref{sec6}).

\subsection{McDonald--Kreitman tables}\label{sec22}
Theorems \ref{theo21} and \ref{theo22} give the
distribution of polymorphic sites and fixations of mutant nucleotides
in a
large population but are not directly applicable to samples from a finite
number of individuals.

Suppose that we have random samples from two species that are related but
not extremely distantly related \cite{rMcDKr}, for example, the
\textit{Drosophila} species \textit{melanogaster} and \textit
{simulans}
\cite{rNatMod,rSSKul,rSSPar} or humans and
chimpanzees \cite{rNielsen}. Specifically, assume that we have a DNA
alignment of $m+n$ sequenced genes from the same genetic locus in two
species, of which $m$ are randomly chosen from one species and $n$ from
the second species. The sample can be viewed as an $(m+n)\times L$ matrix,
where $L$ is the number of nucleotide sites in the alignment and the
matrix elements are chosen from the four letters A, C, G, T, which
stand for the four nucleotides in DNA. If mutation is sufficiently rare so
that it never occurs more than once at the same site, then there will be
at most two nucleotides at any one site in the population and hence at
most two letters in any one column of the matrix.

Given the joint alignment, we say that a site is a \textit{fixed
difference}
if it is monomorphic in the sample within each species but polymorphic in
the two species (i.e., monomorphic within each species, but at
different nucleotides). The site is \textit{polymorphic} if it is polymorphic
within either or both of the two species.

Let $K_s,K_r$ be the number of silent and replacement fixed differences
in the joint alignment, and $V_s,V_r$ the number of (within-species)
polymorphisms at silent and replacement sites. These counts can be
arranged in the $2\times2$ contingency table
%
\begin{equation}\label{equ218}
\matrix{
& \matrix{ D & P} \cr
\matrix{ S\cr R} &
\left[\matrix{K_s & V_s\cr K_r & V_r} \right]},
\end{equation}
where the column headings $D,P$ refer to fixed differences and
polymorphisms and the row headings $S,R$ to silent and replacement
mutations, respectively. The table (\ref{equ218}) is called a McDonald--Kreitman
table \cite{rMcDKr} and also a DPRS table due to the row and column
headings in (\ref{equ218}).

A statistically significant excess in the lower-left corner ($K_r$) of
the table suggests that one or both species has seen significant positive
or directional selection at that gene since the two species diverged. In
contrast, a statistically significant deficit of $K_r$ suggests that,
instead, most new replacement mutations have been subject to strong
negative selection. If the sites are assumed independent, one can apply
traditional $2\times2$ contingency table tests to infer an excess or deficit
of replacement fixed differences \cite{rMcDKr,rDurr}.

We say that a polymorphic site in the joint alignment is a \textit{legacy
polymorphism} if it is the descendent of a site that was polymorphic in
the common ancestral species at the time of divergence. In contrast, a
\textit{new polymorphism} is a polymorphic site that was caused by a mutation
since the time of divergence in one of the two daughter species.

An important difference between the two types of polymorphisms is that
legacy polymorphisms can lead to shared polymorphisms, which are sites
that are polymorphic in both samples. In contrast, if mutations are
sufficiently rare so that they never occur more than once at the same
site, then new polymorphisms can cause sample polymorphisms in only one
daughter species. Another difference is that legacy polymorphisms begin at
mutant population frequencies strictly between 0 and 1, while new
polymorphisms begin at mutant population frequency 0 (or $1/N$ for the
Markov chain). Thus legacy polymorphisms are more likely to have multiple
copies of both the mutant and original nucleotides for relatively recent
divergence times.

An alternative to the $2\times2$ table in (\ref{equ218}), that may have more
accuracy in estimating parameters for recently diverged populations, is
%
\begin{equation}\label{equ219}
\matrix{ & \matrix{ D & O & H}\cr
\matrix{ S\cr R} &
\left[\matrix{K_s & O_s & H_s\cr K_r & O_r & H_r} \right]
},
\end{equation}
where $O_s,O_r$ are numbers of sites in the joint sample that are
polymorphic in only one sample and $H_s,H_r$ the numbers of sites that are
polymorphic in both samples. Thus $V_s=O_s+H_s$ and $V_r=O_r+H_r$. A
disadvantage of the DOHRS table (\ref{equ219}) is that rare events in which
multiple mutations occur at the same site may lead to a pair of new
polymorphisms being misclassified as a legacy polymorphism. On the other
hand, the more usual practice of counting sites that are polymorphic in
both samples as two polymorphic sites may misclassify a single legacy
polymorphism. (See \cite{r67} for other examples of higher-dimensional
McDonald--Kreitman tables.)

\subsection{The distribution of sampling statistics}\label{sec23}
For simplicity, we assume that the effective population sizes $N_e$, the
mutation rates $\theta$ and the selection coefficients $\gamma$ are
the same in
the two daughter populations. It then follows from Theorem \ref{theo31} in
Section \ref{sec33} that, under the assumptions of Theorem \ref{theo21}, the counts
$K_s,V_s,K_r,V_r$ in (\ref{equ218}) and $K_s,O_s,H_s,K_r,O_r,H_r$ in (\ref{equ219}) are
independent Poisson with means depending on $\beta_r=(t,\theta
_r,\gamma
)$ for
replacement sites and $\beta_s=(t,\theta_s)$ at silent sites. If the effective
population sizes for the two daughter species are different, the resulting
formulas in Section \ref{sec3} are easy to modify.

Let $Z_a$ ($1\le a\le4$) be the counts in (for example) (\ref{equ218}) in some
order. Let $m_a=m_a(\beta)=E(Z_a)$ be the corresponding formulas in
Section \ref{sec33} for $\beta=(\beta_r,\beta_s)$. Then the likelihood of (\ref{equ218})
can be
written
%
\begin{equation}\label{equ220}
L(\beta,Z) =\prod_{a=1}^4
\biggl(\exp(-m_a(\beta))
\frac1{Z_a !} m_a(\beta)^{Z_a} \biggr)
\end{equation}
or
%
\begin{equation}\label{equ221}
\log L(\beta,Z) =
C(Z) -\sum_{a=1}^4 m_a(\beta)
+\sum_{a=1}^4 Z_a \log m_a(\beta).
\end{equation}
While in principle the parameters $\beta_r,\beta_s$ can be estimated by
maximizing the log likelihood (\ref{equ221}), normally only a few loci in a few
species are sufficiently polymorphic to allow $\beta$ to be estimated from
data from a single locus. More typically, likelihoods of the form (\ref{equ220})
are combined over many loci to form a single likelihood. The resulting
expression can either be maximized to estimate the model parameters
\cite{rBoyko,rKeight2,rWilmson} or else analyzed by Bayesian methods
\cite{rNatMod,rSSKul,rSSPar,rBaines,rHaley}.

Although the mutation rate per site per generation is normally assumed to
be the same for sites in the same genetic locus, which would suggest
$\theta_r>2\theta_s$, most replacement mutations are either lethal to
their host
or else are severely deleterious. This means that most replacement
mutations are immediately lost on the diffusion time scale. This can show
up in the model as a censored replacement mutation rate with
${\theta_r/\theta_s<1}$ and sometimes ${\theta_r/\theta_s\ll1}$.
In particular, there
is no simple relation between $\theta_s$ and $\theta_r$ in (\ref{equ220}).

Simulations have shown that methods based on (\ref{equ220}) for multilocus data
are relatively robust to violations of basic model assumptions, such as
lack of local independence of polymorphic sites
\cite{rHaley,rBoyko,rDirSel,rZhuLan}.

\subsection{Final comments}\label{sec24}
These results generalize sampling formulas in Sawyer and
Hartl \cite{rSSDH} who made the approximation that the two species
populations are individually at equilibrium. Numerical simulations have
shown that these lead to biased estimates of the population divergence
time \cite{rHaley}, particularly when the divergence time is small, but
that estimates of the mutation and selection parameters are relatively
unbiased. Nonequilibrium results can also be applied to situations
where a
population has been subject in the past to an abrupt change in selection
parameters or population size (see, e.g., \cite{rWilmson}).

The sampling formulas in Section \ref{sec3} for the Poisson means in (\ref{equ218})
and (\ref{equ219}) are more complex than the equilibrium results \cite{rSSDH}, but
likelihoods based on these means can be analyzed numerically as in the
time homogeneous case. Terms involving diffusion transition densities
(Section \ref{sec3}) can be estimated by the Crank--Nicholson method
\cite{rPressEtal,rWilmson}. As in \cite{rSSPar}, Gauss--Legendre
quadrature can be used for integrals over a finite range. Generalizations
for random distributions of selection coefficients within loci have been
handled by Gauss--Hermite quadrature for within-locus normal
variation \cite{rBaines,rSSPar} and Gauss--Laguerre quadrature for
within-locus exponential variation \cite{rHaley}. The behavior of the
models described in this paper on simulated data, and applications to
biological DNA sequence data, will be discussed in future publications.

\section{Sampling formulas for aligned DNA sequences}\label{sec3}
Since the two species in Section \ref{sec22} are assumed to be relatively closely
related on an evolutionary time scale, it is natural to assume that they
have the same mutation and selection rates at each genetic locus and also
the same average generation times. For simplicity, we also assume that the
effective population sizes of the two daughter species are the same. This
means that the scaled mutation rates $\theta_s,\theta_r$ and selection
coefficients $\gamma$ are the same at each genetic locus, and the scaled
time $t$ since the divergence of the two species from their most recent
common ancestral species is also the same [see~(\ref{equ29})]. As mentioned
earlier, the formulas below can be easily modified if these assumptions
are violated.

As in Section \ref{sec22}, assume that we have a DNA alignment of $m+n$ sequenced
genes from the same genetic locus in two species, of which $m$ are
randomly chosen from one species and $n$ from the second species. Each of
the sample statistics $K_s,V_s,K_r,V_r$ in (\ref{equ218}) and
$K_s,O_s,H_s,K_r,O_r,H_r$ in (\ref{equ219}) are counts of sites of two types,
legacy and new polymorphisms, that lead to different formulas for the
expected counts.

\subsection{Legacy polymorphisms}\label{sec31}
It follows from Section \ref{sec71} that the distribution of
polymorphic site frequencies $p_i$ in a single daughter population at time
$t>0$ that are derived from ancestral legacy polymorphic sites is a
Poisson random field with mean density
%
\begin{equation}\label{equ31}
f_L(\beta,y) =\int_0^1 p(t,x,y)\nu(dx),
\end{equation}
where $\beta=\beta_s=(t,\theta_s,0)$ for silent mutations and
$\beta=\beta_r=(t,\theta_r,\gamma)$ for replacement mutations. It
follows similarly from Section \ref{sec8} that the number of population
mutant fixations in a single daughter population by time $t>0$ that are
derived from legacy polymorphic sites is Poisson with mean
\[
G_L(\beta) = \int_0^1 P_x(T_1\le t) \nu(dx).
\]

Assume that a particular legacy polymorphic site has population
frequency $x$ ($0<x<1$) for the mutant nucleotide at time $t=0$. For a
sample of size $n$ from a daughter population at time $t>0$, let $I(x,n)$
be the probability that this site is monomorphic in the sample at the
wild-type (nonmutant) nucleotide, $J(x,n)$ that it is polymorphic in the
sample and $K(x,n)$ the probability that it is monomorphic at the mutant
nucleotide. Then we have the following lemma.
\begin{lem}\label{lemm31} For $I(x,n)$, $J(x,n)$ and $K(x,n)$ defined
above,
%
\begin{eqnarray}\label{equ32}
I(x,n) &=& P_x(T_0\le t) + \int_0^1 p(t,x,y)(1-y)^n m(dy), \nonumber\\
J(x,n) &=& \int_0^1 p(t,x,y)\bigl(1 - y^n - (1-y)^n\bigr) m(dy),
\\
K(x,n) &=& P_x(T_1\le t) + \int_0^1 p(t,x,y)y^n m(dy) \nonumber
\end{eqnarray}
for $T_a=\min\{a\dvtx X_t=a\}$ as in Section \ref{sec2}.
\end{lem}

The first terms in the formulas for $I$ and $K$ are due to population
fixations. The remaining terms are due to sampling from polymorphic sites.
Note that $I+J+K=1$ in Lemma \ref{lemm31} and that $I$, $J$ and $K$ depend
implicitly on $\beta=(t,\theta,\gamma)$ through both $p(t,x,y)$ and $m(dy)$.
\begin{pf*}{Proof of Lemma \ref{lemm31}}
For a legacy polymorphism of initial mutant frequency $x$, condition on
its population frequency $y$ at time $t>0$.
\end{pf*}

Given a sample of size $m+n$ from the two populations at time $t>0$, let
$L_1$, $L_2$, $L_3$ be the random numbers of legacy polymorphism sites that
are fixed differences in the sample ($L_1)$, polymorphic in only one of
the two samples ($L_2$) or polymorphic in both samples ($L_3$). Then we
have the following lemma.
\begin{lem}\label{lemm32}
The random variables $L_1,L_2,L_3$ defined above are
independent Poisson with means
%
\begin{eqnarray}\label{equ33}\quad
E(L_1) &=& C_1(\beta) = \int_0^1 \bigl(I(x,m)K(x,n) +
I(x,n)K(x,m)\bigr)
\nu(dx), \nonumber\\
E(L_2) &=& C_2(\beta) = \int_0^1 \bigl[J(x,m)
\bigl(I(x,n)+K(x,n)\bigr)
\nonumber\\
&&\hspace*{57.2pt}{} + J(x,n)\bigl(I(x,m)+K(x,m)\bigr)\bigr] \nu(dx)
\\
&=& \int_0^1 \bigl(J(x,m) + J(x,n) - 2J(x,m)J(x,n)\bigr)
\nu(dx), \nonumber\\
E(L_3) &=& C_3(\beta) = \int_0^1 J(x,m)J(x,n) \nu(dx). \nonumber
\end{eqnarray}
\end{lem}
\begin{pf}
Consider the random mapping $x\to z=\{1,2,3\}$
corresponding to random sampling at time $t>0$ for these three outcomes.
Then Lemma \ref{lemm32} follows from Lemma \ref{lemm31} and Bartlett's
theorem \cite{rKingman}.
\end{pf}

In particular, it follows from Lemma \ref{lemm32} that the expected number of
polymorphic sites at time $t>0$ due to legacy polymorphisms is
\[
E(L_2+L_3) =
\int_0^1\bigl(J(x,m) + J(x,n) - J(x,m)J(x,n)\bigr) \nu(dx).
\]
The third term in the integrand above corrects for double counting at
shared polymorphisms.

\subsection{Polymorphisms from new mutations}\label{sec32}
It follows from Section \ref{sec72} that the contribution of new\vadjust{\goodbreak}
 polymorphisms
to site polymorphisms at time $t>0$ in a single population is a Poisson
random field with density
%
\begin{equation}\label{equ34}\quad
f_N(\beta,y) = \frac\theta{s(1)} \biggl(s(1) - s(y)
-\int_0^1\bigl(s(1)-s(x)\bigr) p(t,x,y)m(dx) \biggr). \\
\end{equation}
It follows similarly from Section \ref{sec8} that the number of mutant fixations in
a single daughter population due to new polymorphisms is a Poisson random
variable with mean
%
\begin{equation}\label{equ35}
G_N(\beta) = \frac\theta{s(1)}
\int_0^t \widetilde{P}_0(T_1\le u) \,du.
\end{equation}
In particular, the limiting PRF mean density for population site
polymorphic frequencies in a single population is $f(\beta,y) =
f_L(\beta,y)+f_N(\beta,y)$ as in (\ref{equ214}), and the limiting expected
number of
mutant fixed differences is $G(\beta)=G_L(\beta)+G_N(\beta)$ as in (\ref{equ216}).

For a sample of size $n$ from a single population at time $t>0$, let $Z_k$
($0\le k\le n$) be the number of new polymorphic sites that have $k$
mutant nucleotides. Then we have the following lemma.
\begin{lem}\label{lemm33}
The random variables $Z_k$ defined above are independent
Poisson-distributed random variables with means
%
\begin{equation}\label{equ36}
E(Z_k) =F_N(\beta,n,k) =
\int_0^1 f_N(\beta,y) \pmatrix{n\cr k} y^k (1-y)^{n-k} m(dy)
\end{equation}
for $1\le k\le n$ and $f_N(\beta,y)$ in (\ref{equ34}).
\end{lem}
\begin{pf}
Consider the random mapping $y\to K=\{0,1,\ldots,n\}$ defined by binomial
sampling with parameters $y$ and $n$ at each polymorphic site. The range
of the mapping is a Poisson random field by Bartlett's theorem
\cite{rKingman}. Thus the random variables $Z_k$ for $1\le k\le n$ are
independent Poisson with the means in (\ref{equ32}). Note that Bartlett's theorem
applies here even though $Z_0=\infty$ a.s. and $E(Z_0)=\infty$.
\end{pf}

Since the random variables $Z_k$ in Lemma \ref{lemm33} are independent, the
expected number of polymorphic sites at time $t>0$ due to new
polymorphisms is Poisson with mean
%
\begin{equation}\label{equ37}\qquad
E_N(\beta,n) =\sum_{k=1}^{n-1} F_N(\beta,n,k)
=\int_0^1 f_N(\beta,y) \bigl(1 - y^n - (1-y)^n\bigr) m(dy).
\end{equation}
The expected number of sites due to new polymorphisms that are monomorphic
in a sample of size $n$ at the mutant nucleotide is (i) the expected
number of sites that have fixed in the population at the mutant nucleotide
by time $t$ plus (ii) the expected number of sites that are not fixed in
the population but are monomorphic at the mutant nucleotide in a sample,
which is
%
\begin{equation}\label{equ38}
D_N(\beta,n) =G_N(\beta) + F_N(\beta,n,n)
\end{equation}
for $G_N(\beta)$ in (\ref{equ35}).\vadjust{\goodbreak}

\subsection{Sampling formulas in DPRS and DOHRS tables}\label{sec33}

We have now completed the proof of the following result.
\begin{theorem}\label{theo31}
Assume that the two species have the same effective population size $N_e$
and the same scaled parameter values $\beta_r=(t,\theta_r,\gamma)$ and
$\beta_s=(t,\theta_s,0)$. Then the counts $K_s,O_s,H_s,K_r,O_r,H_r$
in the
table (\ref{equ219}) in Section \ref{sec2} are independent and Poisson distributed with
means
%
\begin{eqnarray}\label{equ39}
E(K_s) &=& C_1(\beta_s) + D_N(\beta_s,m) +D_N(\beta_s,n), \nonumber\\
E(O_s) &=& C_2(\beta_s) + E_N(\beta_s,m) + E_N(\beta_s,n), \nonumber
\\
E(H_s) &=& C_3(\beta_s), \nonumber\\[-9pt]\\[-9pt]
E(K_r) &=& C_1(\beta_r) + D_N(\beta_r,m) +D_N(\beta_r,n), \nonumber\\
E(O_r) &=& C_2(\beta_r) + E_N(\beta_r,m) + E_N(\beta_r,n), \nonumber
\\
E(H_r) &=& C_3(\beta_r) \nonumber
\end{eqnarray}
for $C_i(\beta)$ in (\ref{equ33}), $E_N(\beta,n)$ in (\ref{equ37}) and $D_N(\beta)$
in (\ref{equ38}).
\end{theorem}
\begin{corr}\label{corr31}
The counts $K_s,V_s,K_r,V_r$ in the $2\times2$ table (\ref{equ218})
are independent Poisson with means $E(V_r)=E(H_r)+E(O_r)$ and
$E(V_s)=E(H_s)+E(O_s)$ in (\ref{equ39}).
\end{corr}

\section{Diffusion operators and diffusion approximations}\label{sec4}
This section describes the diffusion approximation
\cite{rEwens,rHartlPrimer,rTrotter} of a single Markov chain
$X_{1,a,k}$ ($k\ge0$) or
$X_{2,b,r,k}$ ($k\ge r\ge1$) in Section \ref{sec2}, which are assumed to have the
same transition function $p_N(i,j)$ in (\ref{equ21}). A major purpose of this
section is to show that the transition density $p(t,x,y)$ of the limiting
diffusion process and its first partial derivative
$(\partial/\partial s(x))p(t,x,y)$ are smooth for $t>0$ and $0\le
x,y\le1$. (The latter will be used in Section \ref{sec82}.)

Each of the Markov chains $X_{1,a,k}, X_{2,b,r,k}$ will be
approximated by
a diffusion process $X_t$ for continuous $t$ scaled in terms of $N^2$
steps of the Markov chain, so that $t\sim k/N^2$. The limiting diffusion
process $X_t$ is determined by the differential operator
%
\begin{equation}\label{equ41}
L_x =x(1-x)\,\frac{d^2}{dx^2} +\gamma x(1-x)\,\frac{d}{dx},
\end{equation}
where $\gamma=\lim_{N\to\infty} N\sigma_N$ for $\sigma_N$
in (\ref{equ21}). We write the
operator $L_x$ in Feller form $L_x=(d/m(dx))
(d/s(dx))$
where
%
\begin{equation}\label{equ42}
s(x)=\frac{1-e^{-\gamma x}}{\gamma} \quad\mbox{and}\quad
m(dx)=\frac{e^{\gamma x}}{x(1-x)}\,dx
\end{equation}
are called the \textit{scale} and \textit{speed measure} of $L_x$,
respectively
\cite{rEwens,r38Feller,rItoMcK,rKarlin,rSSFatou}. At silent sites, and
in general if $\gamma=0$, the scale and speed measure are $s(x)=x$ and
$m(dx)=dx/(x(1-x))$.

\subsection{Green's functions and transition densities}\label{sec41}
Define
%
\begin{equation}\label{equ43}
g(x,y) =
\frac{(s(1)-s(x\vee y))(s(x\wedge y)-s(0))}
{s(1)-s(0)},
\end{equation}
where $x\vee y=\max\{x,y\}$ and $x\wedge y=\min\{x,y\}$. Set
%
\begin{equation}\label{equ44}
B_{01} = \{f\in C[0,1]\dvtx f(0)=f(1)=0\},
\end{equation}
where $C[0,1]$ is the class of continuous functions on $0\le x\le1$.
Then, for any ${f\in C[0,1]}$, $h(x)=\int_0^1 g(x,y)f(y)m(dy)$ is the
unique solution $h(x)\in C^2(0,1)\cap B_{01}$ such that $L_xh(x)=-f(x)$
for $0<x<1$. (For $s(x),m(dx)$ in (\ref{equ42}), this does \textit{not} imply
$h\in
C^2[0,1]$, nor even that $h\in C^1[0,1]$.)

Since $s(x)$ is increasing, $g(x,y) \le\min\{g(x,x), g(y,y)\}$
by (\ref{equ43}). Hence by (\ref{equ43}) and (\ref{equ42})
%
\begin{equation}\label{equ45}
k(x) = \int_0^1 g(x,y)^2m(dy) \le g(x,x)\int_0^1 g(y,y)m(dy) <
\infty
\end{equation}
and $\int_0^1 k(x)m(dx) \le(\int_0^1 g(x,x)m(dx))^{ 2}<\infty$. Thus
\[
\int_0^1\int_0^1 g(x,y)^2m(dx)m(dy) <\infty
\]
and $g(x,y)$ is a Hilbert--Schmidt kernel on $L^2(I,m)$ \cite{rRieszNagy}.
This implies that there exists a complete orthonormal system of functions
$\alpha_n(x)$ in $L^2(I,m)$ such that
%
\begin{equation}\label{equ46}
\int_0^1 g(x,y)\alpha_n(y)m(dy) = \beta_n \alpha_n(x),\qquad
\int_0^1 \alpha_n(y)^2m(dy)=1,
\end{equation}
where $\sum_{n=1}^\infty\beta_n^2<\infty$ and $\alpha_n(0)=\alpha
_n(1)=0$. By the
same arguments as after~(\ref{equ44}), $\beta_n\ne0$ and the integral equation
in (\ref{equ46}) is equivalent to
%
\begin{equation}\label{equ47}
L_x\alpha_n(x) = -\lambda_n\alpha_n(x),\qquad
\lambda_n=1/\beta_n,\qquad
\alpha_n\in C^2(0,1)\cap B_{01}.
\end{equation}
Since $\alpha_n(0)=\alpha_n(1)=0$, $\lambda_n>0$ by (\ref{equ41}) and $\sum
_{n=1}^\infty
1/\lambda_n^2<\infty$. In particular, we can assume
$0<\lambda_1\le\lambda_n\uparrow\infty$.

By Cauchy's inequality in (\ref{equ46}) and by (\ref{equ45}) and (\ref{equ43})
\[
|\alpha_n(x)| \le\lambda_n\sqrt{k(x)}
\le C_1\lambda_n\sqrt{x(1-x)} \le C_1\lambda_n.
\]
Since $g(x,y)\le g(x,x)$ and $\int_0^1 \sqrt{y(1-y)}m(dy)<\infty$, (\ref{equ46})
implies
%
\begin{equation}\label{equ48}
|\alpha_n(x)| \le C_2\lambda_n^2 g(x,x)
\le C_3\lambda_n^2 x(1-x).
\end{equation}
Since $g(x,y)/(x(1-x))$ is bounded and continuous by (\ref{equ43}), it follows
from (\ref{equ46}) and (\ref{equ48}) that $\alpha_n(x)/(x(1-x))$ extend to continuous
functions on $[0,1]$. If $\gamma=0$, the functions $\alpha_n(x)$ are
polynomials
related to Jacobi polynomials \cite{rKimura}. If $\gamma\ne0$, they are
entire functions that are never polynomials.

By Mercer's theorem \cite{rRieszNagy},
%
\begin{equation}\label{equ49}
g(x,y)=\sum_{n=1}^{\infty}\frac{\alpha_n(x)\alpha_n(y)}{\lambda
_n}
\end{equation}
converges absolutely and uniformly for $0\le x,y\le1$. The series
%
\begin{equation}\label{equ410}
p(t,x,y) = \sum_{n=1}^{\infty}e^{-\lambda_nt}\alpha_n(x)\alpha
_n(y)
\end{equation}
converges uniformly for $0\le x,y\le1$ and $t\ge a>0$ since
$\sum_{n=1}^\infty1/\lambda_n^2<\infty$. Thus
%
\begin{equation}\label{equ411}
\int_a^\infty p(u,x,y)\,du
= \sum_{n=1}^{\infty}e^{-\lambda_n a}
\frac{\alpha_n(x)\alpha_n(y)}{\lambda_n}
\end{equation}
converges absolutely and uniformly for $a\ge0$. In particular
%
\begin{equation}\label{equ412}
g(x,y)=\int_0^\infty p(t,x,y) \,dt
\end{equation}
with uniform convergence for $0\le x,y\le1$.

It follows from (\ref{equ410}) and (\ref{equ47}) that
%
\begin{eqnarray}\label{equ413}
p(t+s,x,y) &=& \int_0^1 p(t,x,z)p(s,z,y)m(dz)\quad
\mbox{and}\nonumber\\[-8pt]\\[-8pt]
(\partial/\partial t)p(t,x,y)
&=& L_x p(t,x,y),\qquad t>0, 0<x,y<1.\nonumber
\end{eqnarray}
Choose\vspace*{1pt} $f\in B_{01}$ with $f(x)=0$ for $0\le x\le c$ and $1-c\le x\le1$
for some $c>0$. Let\vspace*{1pt} $u(t,x)=\int_0^1 p(t,x,y)h(y)m(dy)$ for
$h(x)=\int_0^1g(x,y)f(y)m(dy)$. Then by (\ref{equ411}) and (\ref{equ412})
\[
u(t,x) = \sum_{n=1}^\infty
\frac{e^{-\lambda_n t}}{\lambda_n} \alpha_n(x) c_n \qquad\mbox
{for }
c_n=\int_0^1 f(y)\alpha_n(y) m(dy),
\]
where the series converges uniformly for $0\le x\le1$ and $0\le
t<\infty$
by (\ref{equ411}). Thus $u(t,x)\in C([0,\infty)\times[0,1])$ with
\begin{eqnarray*}
(\partial/\partial t)u(t,x) &=& L_x u(t,x),\qquad t>0, 0<x<1, \\
u(t,0)&=&u(t,1)=0,\qquad u(0,x)=h(x).
\end{eqnarray*}
It follows from maximum principles for parabolic partial differential
equations \cite{rProttWein} that
%
\begin{equation}\label{equ414}
p(t,x,y)\ge0,\qquad \int_0^1 p(t,x,z)m(dz)\le1
\end{equation}
for $0<x,y<1$. If $u(x,t)=Q_th(x)=\int p(t,x,y)h(y)m(dy)$, then
%
\begin{equation}\label{equ415}
\lim_{t\to0} Q_tf(x) = f(x)\vadjust{\goodbreak}
\end{equation}
uniformly for $0\le x\le1$ for a dense set of $f\in C^2(0,1)\cap B_{01}$,
and hence for all $f\in B_{01}$ by standard arguments.

\subsection{Diffusion processes and semigroups}\label{sec42}
It follows from the relations (\ref{equ413})--(\ref{equ415}) that there exists a diffusion
process $X_t$ with continuous sample paths with $0\le X_t\le1$ such that
%
\begin{equation}\label{equ416}
Q_tf(x) = E_x( f(X_t)) = \int_0^1 p(t,x,y)f(y)m(dy)
\end{equation}
for $f\in C[0,1]$ with $f(0)=f(1)=0$ \cite
{rDynkin,rEthierKurtz,r37Feller,rItoMcK,rSSFatou}.
The process $X_t$ satisfies $0<X_t<1$ up
to the time that it is trapped at one of the endpoints $0,1$, which happens
eventually with probability one. The relation
%
\begin{equation}\label{equ417}
|Q_tf(x)| \le C_4 e^{-\lambda_1t} \|f\|,\qquad
\|f\|={\sup_{0\le y\le1}} |f(y)|
\end{equation}
from (\ref{equ410}) and (\ref{equ48}) gives the rate at which $X_t$ is trapped at the
endpoints. The exit point for $X_t$ is given by the scale function
%
\begin{equation}\label{equ418}
P_x(T_1<T_0) =\frac{s(x)-s(0)}{s(1)-s(0)}
=\frac{s(x)}{s(1)},
\end{equation}
where $T_a=\min\{s\dvtx X_s=a\}$ \cite{r37Feller,rItoMcK}.

It follows from (\ref{equ48})--(\ref{equ415}) that $Q_t\dvtx B_{01}\to B_{01}$ for all $t>0$
for $B_{01}$ in (\ref{equ44}) and that $\{Q_t\}$ is a strongly continuous
semigroup of linear operators on $B_{01}$.

The \textit{infinitesimal generator} \cite{rDunfSchw,rRieszNagy,rTrotter}
of a semigroup of linear operators $Q_t$ on a Banach space $B$ is the
linear operator $A$ defined by $Ah=f$ on the linear subspace
\[
\mathcal{D}(A)=\Bigl\{h\in B\dvtx \lim_{t\to0}\|(1/t)(Q_th-h)-f\| =0
\mbox{ for some $f\in B$}\Bigr\},
\]
where $\|f\|$ is the norm in the Banach space. If $Q_t$ is strongly
continuous, then $\mathcal{D}(A)$ is dense in $B$ and
%
\begin{equation}\label{equ419}
\|Q_tf\|\le Me^{Kt}\|f\| \qquad\mbox{all }f\in B
\end{equation}
for some real $K$. If $K<0$ as in (\ref{equ417}), $\mathcal{D}(A)$ is the
range of the
resolvent operator $R_0 f=\int_0^\infty Q_tf\,dt$ on $B$ with
$-AR_0=I$.
This implies $Ah=-f$ if $h=R_0f$ \cite{rDunfSchw,rRieszNagy,rTrotter}.
By (\ref{equ412})
%
\begin{equation}\label{equ420}
R_0f(x) = \int_0^\infty Q_tf(x) \,dt = \int_0^1 g(x,y)f(y)m(dy)
\end{equation}
is the Green's operator defined by $g(x,y)$.

A \textit{core} of a strongly-continuous semigroup $Q_t$ with $K<0$
is a
subset ${\mathcal{C}\subseteq\mathcal{D}(A)}$ such that
$B_c=A(\mathcal{C})$ is dense in $B$. Since
$R_0$ is one--one on $B$, this is equivalent to specifying a dense subset
$B_c\subseteq B$ and setting $\mathcal{C}=R_0(B_c)$.

\subsection{Diffusion approximations and Trotter's theorem}\label{sec43}
Let $X_k^N$ be the Markov chain defined by the Moran model (\ref{equ21}). Define
$Y^N_j=X^N_j/N$, so that $0\le Y^N_j\le1$. It follows from standard
arguments and (\ref{equ21}) that
\begin{theorem}\label{theo41}
Let $i_N$ be integers such that $0\le i_N\le N$ and such that
$x_N=i_N/N\to x$ for some $x$, $0\le x\le1$. Then for any $\delta>0$
%
\begin{eqnarray}\label{equ421}
\lim_{N\to\infty}N^2E_{i_N} (Y^N_1 - x_N)
&=& \gamma x(1-x),\nonumber\\
\lim_{N\to\infty}N^2E_{i_N} \bigl((Y^N_1 - x_N
)^2 \bigr)
&=& 2x(1-x), \\
\lim_{N\to\infty}N^2E_{i_N}
(|Y^N_1 - x_N|^{2+\delta} ) &=& 0.\nonumber
\end{eqnarray}
\end{theorem}

Since the functions on the right-hand side of (\ref{equ421}) are continuous,
convergence in (\ref{equ421}) is equivalent to uniform convergence in $x$ for
${i_N=[Nx]}$. By Taylor's theorem
%
\begin{eqnarray}\label{equ422}
&&\lim_{N\to\infty}N^2 E_{i_N} \bigl(h(Y^N_1)-h(x_N) \bigr)
\nonumber\\[-8pt]\\[-8pt]
&&\qquad =L_xh(x) =x(1-x)h''(x) + \gamma x(1-x)h'(x)
\nonumber
\end{eqnarray}
uniformly for $0\le x\le1$ for any $h \in C^2[0,1]$. Then by Trotter's
theorem \cite{rTrotter}:
\begin{theorem}\label{theo42}
For $Y^N_j$ as above, $i_N=[Nx]$, and the diffusion
process $X_t$ in (\ref{equ416}),
%
\begin{equation}\label{equ423}
\lim_{N\to\infty}E_{i_N} \bigl(f\bigl(Y^N_{[N^2t]}\bigr) \bigr)
= E_x(f(X_t)) = Q_tf(x)
\end{equation}
uniformly for $0\le x\le1$ for any $f\in C[0,1]$ with $f(0)=f(1)=0$. The
convergence is also uniform for $0\le t\le T$ for any $T>0$.
\end{theorem}

We cannot apply (\ref{equ422}) for $h \in C^2[0,1]$ directly for Trotter's
theorem, since in this case there exist $h\in\mathcal{D}(A)$ with
$h\notin
C^1[0,1]$, let alone $C^2[0,1]$. However, it is sufficient to
verify (\ref{equ422}) for all $h$ in a core for $A$ \cite{rTrotter}. If
$\mathcal{C}=R_0(B_c)$ where $B_c$ is the set of all function $f\in
B_{01}$ such
that $f(x)=0$ for $0\le x\le a$ and $1-a\le x\le1$ for some $a>0$, then
$\mathcal{C}$ is such a core.

The result (\ref{equ423}) also holds for $f\in C[0,1]$ without the conditions
$f(0)=f(1)=0$ with an appropriate modification of the definition
of $Q_tf(x)$. See Corollary \ref{corr51} in Section \ref{sec5} below.

Thus, after suitable rescaling, the Markov chains $\{X_{1,a,k},
X_{2,b,r,k}\}$ in Section \ref{sec2} converge in distribution to diffusion
processes $\{X_t\}$ with infinitesimal generator (\ref{equ25}) and scale and
speed measure (\ref{equ27}) in Section \ref{sec2}.

\section{Exit probabilities for Markov chains}\label{sec5}
Let $X_k^N$ be the Moran-model Markov chain defined by (\ref{equ21}). Recall that
$N\sigma_N\to\gamma$ as $N\to\infty$ by (\ref{equ29}). Then we have
the following lemma.
\begin{lem}\label{lemm51}
Let $i,m$ be integers such that $1\le i\le m\le N$. Then
%
\begin{equation}\label{equ51}
P_i(T^N_m < T^N_0) =
\frac{ 1 - (1+\sigma_N)^{-i}}{1 - (1+\sigma_N)^{-m}}
\end{equation}
with the right-hand side replaced by $i/m$ if $\sigma_N=0$.
\end{lem}
\begin{pf}
By (\ref{equ21}), $p_N(i,i+1)/p_N(i,i-1) = 1+\sigma_N$ for $0<i<N$. Since we can
ignore ``wait states'' with $X_{k+1}^N=X_k^N$, (\ref{equ51}) follows from the
classical Gambler's Ruin problem (see, e.g., \cite{rKarlin}, pages 50, 92--94,
and Moran \cite{rMoran}). (See Lemma \ref{lemm82} for a second proof.)
\end{pf}

As a consequence of Lemma \ref{lemm51}:
\begin{lem}\label{lemm52}
Let $i_N$ be integers such that $0\le i_N\le N$ and $i_N/N\to x$ for
some $x$, $0\le x\le1$. Then
%
\begin{equation}\label{equ52}
\lim_{N\to\infty}P_{i_N} (T^N_N<T^N_0 )
=P_x (T_1<T_0 ) =\frac{s(x)-s(0)}{s(1)-s(0)} =
\frac{s(x)}{s(1)}
\end{equation}
for $T_0,T_1$ in (\ref{equ418}) and $s(x)$ in (\ref{equ42}). If $i_N=i_N(x)=[Nx]$, the
convergence in (\ref{equ52}) is uniform in $x$ for $0\le x\le1$.
\end{lem}
\begin{pf}
If $\sigma_N\ne0$ and $\gamma\ne0$, it follows from (\ref{equ51}) that
%
\begin{equation}\label{equ53}
P_{i_N}(T^N_N < T^N_0) =
\frac{1 - (1+{N\sigma_N}/{N} )^{-i_N}}
{1 - (1+{N\sigma_N}/{N} )^{-N}}
\to \frac{1-e^{-\gamma x}}{1-e^{-\gamma}}
\end{equation}
as $N\to\infty$. The proof is similar if $\gamma=0$.
\end{pf}

Lemma \ref{lemm52} can be used to extend Theorem \ref{theo42} to all $f\in C[0,1]$:
\begin{corr}\label{corr51} Assume $f\in C[0,1]$. Let $X^N_k$ and set $Y^N_k=X^N_k/N$ as
in Theorem \ref{theo41}. Set ${i_N=[Nx]}$. Then
%
\begin{equation}\label{equ54}
\lim_{N\to\infty}E_{i_N} \bigl(f
\bigl(Y^N_{[N^2t]}
\bigr) \bigr)
= E_x(f(X_t)) = \overline{Q}_tf(x)
\end{equation}
uniformly for $0\le x\le1$, where
%
\begin{equation}\label{equ55}
\overline{Q}_tf(x) = f(0)P_x(T_0\le t) + \int_0^1 p(t,x,y)f(y)m(dy)
+ f(1)P_x(T_1\le t).\hspace*{-32pt}
\end{equation}
\end{corr}
\begin{pf}
Any $f\in C[0,1]$ can be written
\[
f(x) = g(x) + f(0)\frac{s(1)-s(x)}{s(1)-s(0)}
+ f(1)\frac{s(x)-s(0)}{s(1)-s(0)},
\]
where $g(0)=g(1)=0$, so that $g\in B_{01}$. Since (\ref{equ54}) holds for $g\in
B_{01}$ by Theorem~\ref{theo42}, it only remains to prove (\ref{equ54}) for $f(x)=s(x)$.
However,
%
\begin{equation}\label{equ56}
h_N(i_N) =P_{i_N}(T^N_N < T^N_0) =E_{i_N}(h(X^N_1))
=E_{i_N}\bigl(h_N\bigl(X^N_{[N^2t]}\bigr)\bigr)
\end{equation}
for all $t>0$. Lemma \ref{lemm52} applied to both sides of (\ref{equ56}) implies (\ref{equ54}) for
$f(x)=s(x)$.
\end{pf}

A stronger result than Lemma \ref{lemm52} is the ``local limit theorem.''
\begin{lem}\label{lemm53} Let $i_N$ be integers with $1\le i_N\le N$
and set
$x_N=i_N/N$. Then
%
\begin{equation}\label{equ57}
\lim_{N\to\infty}
\frac{P_{i_N}(T^N_N < T^N_0)}{s(x_N)/s(1)} =1.
\end{equation}
Similarly, if $0\le i_N\le N-1$ and $x_N=i_N/N$, then
%
\begin{equation}\label{equ58}
\lim_{N\to\infty}
\frac{1 - P_{i_N}(T^N_N < T^N_0)}{(s(1)-s(x_N))/s(1)} =1.
\end{equation}
\end{lem}

The most important cases of Lemma \ref{lemm53} are when $x_N\to0$ or $x_N\to1$.
If $i_N=i_N(x)=\min\{[Nx]+1,1\}$, then (\ref{equ57}) holds uniformly in $x$ for
$0\le x\le1$. Similarly, (\ref{equ58}) holds uniformly in $x$ if
$i_N=i_N(x)=[Nx]$.
\begin{pf*}{Proof of Lemma \ref{lemm53}}
The ratio in (\ref{equ57}) can be written
%
\begin{eqnarray}\label{equ59}
&&\biggl(\frac{1 - (1+\sigma_N)^{-i_N}}{1-e^{-\gamma x_N}}\biggr)
\biggl(\frac{1-e^{-\gamma}}{1 - (1+\sigma_N)^{-N}}\biggr) \nonumber\\[-8pt]\\[-8pt]
&&\qquad=
\biggl(\frac{\int_0^{x_N} e^{-yN\log(1+\sigma_N)} \,dy}
{\int_0^{x_N} e^{-y\gamma} \,dy}\biggr)
\biggl(\frac{\int_0^1 e^{-y\gamma}\,dy}{\int_0^1 e^{-yN\log(1+\sigma_N)}\,dy}\biggr),
\nonumber
\end{eqnarray}
where the first line of (\ref{equ59}) holds for $\gamma\ne0$ and the second
line for
all $\gamma$. Since $N\sigma_N\to\gamma$, then $N\log(1+\sigma
_N)\to\gamma$ and the
ratios converge uniformly in $x_N$. Thus (\ref{equ57}) follows from (\ref{equ59}).
Similarly,
\[
P_x(T_0<T_1) =\frac{s(1)-s(x)}{s(1)-s(0)}
=\frac{e^{-\gamma x}-e^{-\gamma}}{1-e^{-\gamma}}
=\frac{e^{\gamma(1-x)}-1}{e^\gamma-1}
\]
with a corresponding relation for $h_N(i)$. Then the ratio in (\ref{equ58}) can be
written
\[
\biggl(\frac{\int_0^{1-x_N} e^{yN\log(1+\sigma_N)} \,dy}
{\int_0^{1-x_N} e^{y\gamma} \,dy}\biggr)
\biggl(\frac{\int_0^1 e^{y\gamma}\,dy}{\int_0^1 e^{yN\log(1+\sigma_N)}\,dy}\biggr)
\]
with a similar conclusion.
\end{pf*}

\section{Dual Markov chains and dual diffusion processes}\label{sec6}
Define $h_N(i)=P_i(T^N_N<T^N_0)$ as in (\ref{equ51}) where, as before, $P_i$
means conditional on ${X^N_0=i}$. For $p_N(i,j)=P_i(X^N_1=j)$
in (\ref{equ21}) and $1\le i\le N$, define
%
\begin{equation}\label{equ61}
q_N(i,j) =P_i (X_1^N=j\mid T^N_1<T^N_0 )
=\frac1{h_N(i)}p_N(i,j)h_N(j).
\end{equation}
Since $\sum_{j=1}^N q_N(i,j)=1$ for $1\le i\le N$, $q_N(i,j)$ defines a
Markov chain $\{\widetilde{X}^N_k\}$ on $S_N=\{1,2,\ldots,N\}$ that
never attains
$\widetilde{X}^N_k=0$ and has $N$ as an absorbing boundary. Similarly
\[
E_i (\Phi(\widetilde{X}^N_1,\ldots,\widetilde{X}^N_k) )
=E_i \bigl(\Phi(X^N_1,\ldots,X^N_k)\mid T_1<T_0 \bigr)
\]
for all functions $\Phi(j_1,\ldots,j_k)$ on $(S_N)^k$. The chain
$\widetilde{X}_k^N$
[or $q_N(i,j)$] can be called an \textit{$h$-process} of
$X_k^N$ \cite{rKemenySnell}.

\subsection{The limiting dual diffusion process}\label{sec61}
We use the formula for $h_N(i)$ in Lemma \ref{lemm51} to find a diffusion approximation
for $\{\widetilde{X}^N_k\}$. Define $\widetilde{Y}^N_j=\widetilde
{X}^N_j/N$, so that
$0\le
\widetilde{Y}^N_j\le1$. The analog of Theorem \ref{theo41} is:
\begin{theorem}\label{theo61}
Let $i_N$ be integers such that $1\le i_N\le N$ and
$x_N=i_N/N\to x$ where $0\le x\le1$. Then, for any $\delta>0$,
%
\begin{eqnarray}\label{equ62}
\lim_{N\to\infty}N^2\widetilde{E}_{i_N}
(\widetilde{Y}^N_1 - x_N\mid T^N_N<T^N_0 )
&=& b(x) = \gamma x(1-x)\frac{1+e^{-\gamma x}}{1-e^{-\gamma x}},\nonumber\\
\lim_{N\to\infty}N^2\widetilde{E}_{i_N}
\bigl((\widetilde{Y}^N_1 - x_N)^2\mid
T^N_N<T^N_0 \bigr)
&=& 2x(1-x),\\
\hspace*{28pt}\lim_{N\to\infty}N^2\widetilde{E}_{i_N}
(|\widetilde{Y}^N_1 - x_N|^{2+\delta}
\mid T^N_N<T^N_0 ) &=& 0,\nonumber
\end{eqnarray}
where $b(x)=2$ if $x=0$. If $i_N = i_N(x) = \min\{[Nx]+1,1\}$, the
convergence is uniform in $x$.
\end{theorem}
\begin{pf}
By (\ref{equ21}) and Lemma \ref{lemm51}, writing
$j=j_N=i_N$ for ease of notation,
\begin{eqnarray*}
&&N^2 \widetilde{E}_j (Y^N_1 - x_N\mid T^N_N<T^N_0 ) \\
&&\qquad=\frac{N j/N(1-j/N)}
{1+\sigma_N j/N}
\biggl[ \frac{(1+\sigma_N)(1-(1+\sigma_N)^{-j-1}) -
(1-(1+\sigma_N)^{-j+1})}{1-(1+\sigma_N)^{-j}}
\biggr]\\
&&\qquad=\frac{j/N(1-j/N)N}
{1+\sigma_N j/N}
\biggl[ \frac{\sigma_N + \sigma_N(1+\sigma_N)^{-j}}
{1-(1+\sigma_N)^{-j}} \biggr]\\
&&\qquad \to b(x) = \gamma x(1-x)\frac{1+e^{-\gamma
x}}{1-e^{-\gamma x}}.
\end{eqnarray*}
Similarly
\begin{eqnarray*}
&& N^2 \widetilde{E}_j \bigl((Y^N_1 - x_N)^2\mid
T^N_N<T^N_0 \bigr) \\
&&\qquad =\frac{j/N(1-j/N)}
{1+\sigma_N j/N}
\biggl[ \frac{(1+\sigma_N)(1-(1+\sigma_N)^{-j-1}) +
(1-(1+\sigma_N)^{-j+1})}{1-(1+\sigma_N)^{-j}}
\biggr]\\
&&\qquad =\frac{j/N(1-j/N)}
{1+\sigma_N j/N}
\biggl[ \frac{2+\sigma_N - (2+\sigma_N)(1+\sigma_N)^{-j}}
{1-(1+\sigma_N)^{-j}} \biggr]\\
&&\qquad \to 2x(1-x)
\end{eqnarray*}
and
\begin{eqnarray*}
&&N^2 \widetilde{E}_j (|Y^N_1 - x_N|^{2+\delta}\mid
| T^N_N<T^N_0
) \\
&&\qquad =\frac1{N^\delta}
N^2 E_{j/N} \bigl((Y^N_1 - x_N)^2\mid
T^N_N<T^N_0 \bigr)
\to 0.
\end{eqnarray*}
\upqed\end{pf}

As in Section \ref{sec4}, Theorem \ref{theo61} implies, by Taylor's theorem, that
%
\begin{eqnarray}\label{equ63}
&&\lim_{N\to\infty}
N^2\widetilde{E}_j \bigl(h(\widetilde{Y}^N_1)-h(x_N)\mid
T^N_N<T^N_0 \bigr) \nonumber\\[-8pt]\\[-8pt]
&&\qquad
=\widetilde{L}_xh(x) =x(1-x)h''(x) + b(x) h'(x)
\nonumber
\end{eqnarray}
uniformly for $0\le x\le1$ for any $h\in C^2[0,1]$ for $b(x)$ defined
in (\ref{equ62}).

As suggested by $q_N(i,j)=h_N(i)^{-1}p_N(i,h)h_N(j)$ in (\ref{equ61}), the
operator $\widetilde{L}_x$ in (\ref{equ63}) satisfies
\[
\widetilde{L}_xh(x) = \frac1{s(x)}L_x(sh)(x),\qquad 0<x<1.
\]
The operator $\widetilde{L}$ can be written in Feller form
as
\begin{equation}\label{equ64}
\widetilde{L}_xh(x)=\frac{d}{\widetilde{m}(dx)}\frac{d}{\widetilde
{s}(dx)}h(x)
\end{equation}
for scale and speed measure
%
\begin{equation}\label{equ65}
\widetilde{s}(x)=-\frac{1}{s(x)} \quad\mbox{and}\quad
\widetilde{m}(dx)=s(x)^2m(dx)
\end{equation}
for $s(x)$ and $m(dx)$ in Section \ref{sec4}.

Let $\widetilde{X}_t$ be the diffusion process in $(0,1)$ generated
by $\widetilde{L}_x$.
Since
\[
\lim_{x\to0} \widetilde{s}(x)=-\infty\quad\mbox{and}\quad
\int_0^{1/2}|\widetilde{s}(x)|\widetilde{m}(dx)<\infty
\]
the boundary point 0 is an \textit{entrance boundary} for $\widetilde
{L}_x$ or $\widetilde{X}_t$
\cite{r37Feller,rItoMcK,rKarlin}. Since 0\vspace*{1pt} is an entrance
boundary,\vspace*{1pt}
$\widetilde{P}_x(\widetilde{T}_0<\infty)=0$ for $x>0$ where
$\widetilde{T}_a =
\min\{{t\ge0}\dvtx{\widetilde{X}_t=a}\}$ (i.e., 0 is inaccessible
for $\widetilde{X}
_t$), but
\[
\mathop{\lim_{x>0}}_{x\to0} \widetilde{E}_x(f(\widetilde
{X}_t)) =
\widetilde{E}_0(f(\widetilde{X}_t))
\]
exists for all $f\in C[0,1]$ with a probability distribution
$\widetilde{P}_0(\widetilde{X}_t\in dx)$ on $(0,1]$ for $t>0$. Thus
$\widetilde{P}_0(\widetilde{X}_t>0)=1$
for any
$t>0$ and, given $\widetilde{X}_0=0$, $\widetilde{X}_t$ leaves 0
immediately \cite{rItoMcK}.

The Green function for (\ref{equ64}) and (\ref{equ65}) is
%
\begin{eqnarray}\label{equ66}
\widetilde{g}(x,y) &=&
\lim_{\varepsilon\to0}
\frac{(\widetilde{s}(1)-\widetilde{s}(x\vee y))
(\widetilde{s}(x\wedge y)-\widetilde{s}(\varepsilon
))}
{\widetilde{s}(1)-\widetilde{s}(\varepsilon)}\nonumber \\
&=& \widetilde{s}(1)-\widetilde{s}(x\vee y)
=\frac1{s(x\vee y)} -\frac1{s(1)} \\
&=& \frac{(s(1)-s(x\vee y)(s(x\wedge y)-s(0))}
{s(1)s(x\vee y) s(x\wedge y)}
=\frac{g(x,y)}{s(x)s(y)} \nonumber
\end{eqnarray}
by (\ref{equ43}). Since $\iint\widetilde{g}(x,y)^2\widetilde
{m}(dx)\widetilde{m}(dy)
= \iint g(x,y)^2m(dx)m(dy) < \infty$, the kernel $\widetilde
{g}(x,y)$ is
Hilbert--Schmidt with respect to $\widetilde{m}(dx)$. The
eigenfunction equation
\[
\int_0^1\widetilde{g}(x,y)\widetilde{\alpha}(y)\widetilde{m}(dy)
=
\frac1{s(x)}\int_0^1 g(x,y)s(y)\widetilde{\alpha}(y)m(dy)
=\beta\widetilde{\alpha}(x)
\]
is the same as (\ref{equ46}) for $\alpha(x)=s(x)\widetilde{\alpha}(x)$, so
that we can take
$\widetilde{\alpha}_n(x)=\alpha_n(x)/s(x)$ with the same
eigenvalues $\lambda_n$. Define
%
\begin{equation}\label{equ67}
B_1 = \{f\in C[0,1]\dvtx f(1)=0\}
\end{equation}
and
%
\begin{eqnarray}\label{equ68}
q(t,x,y) &=& \sum_{n=1}^\infty e^{-\lambda_n t}\widetilde{\alpha
}_n(x)\widetilde{\alpha}_n(y)
=\frac{p(t,x,y)}{s(x)s(y)}, \nonumber\\
\widetilde{Q}_tf(x)
&=& \widetilde{E}_x(f(\widetilde{X}_t))
=\int_0^1 q(t,x,y)f(y)\widetilde{m}(dy), \\
\widetilde{g}(x,y) &=& \int_0^\infty q(t,x,y)\,dt
=\frac{g(x,y)}{s(x)s(y)} .\nonumber
\end{eqnarray}
For $f\in B_1$, by (\ref{equ68}) and (\ref{equ48}),
%
\begin{equation}\label{equ69}
|\widetilde{Q}_tf(x)| \le C_5 e^{-\lambda_1t} \|f\|,\qquad
\|f\|={\sup_{0\le y\le1}} |f(y)|.
\end{equation}
Since $\widetilde{\alpha}_n(x) = \alpha_n(x)/s(x) \in B_1$ by the
discussion in Section \ref{sec4},
the operators
%
\begin{equation}\label{equ610}
\widetilde{Q}_tf(x)
= \frac1{s(x)} \int_0^1 p(t,x,y)s(y)f(y)m(dy)
= \frac1{s(x)} Q_t(sf)(x)
\end{equation}
preserve $B_1$. The principal result of this section is presented in Section
\ref{sec62}.

\subsection{Trotter's theorem for the dual process}\label{sec62}
Let $\widetilde{X}^N_k$ and $q_N(i,j)$ be for the dual Markov chain in
(\ref{equ61}) and set $\widetilde{Y}^N_k=\widetilde{X}^N_k/N$.

As in Section \ref{sec4}, the relation (\ref{equ63}) defining $\widetilde{L}h$ with uniform
convergence for $h\in C^2[0,1]$ does not cover all
$h\in\mathcal{D}(\widetilde{A})=\widetilde{R}_0(B_1)$ since
$\mathcal{D}(\widetilde{A})$ contains functions $h\notin
C^1[0,1]$. However, if $B_c$ is the set of all $f\in B_1$ such that
$f(x)=b$ for $0\le x\le a$ and $f(x)=0$ for $1-a\le x\le1$ for constants
$a,b$ with $a>0$, then (\ref{equ63}) holds for all $h$ in the core
$\mathcal{C}=\widetilde{R}(B_c)$. Then by Trotter's theorem
\cite{rTrotter}:
\begin{theorem}\label{theo62}
Let $i_N$ be integers such that $1\le i_N\le N$ and
$i_N/N\to x$ for some $x$, $0\le x\le1$. Then
%
\begin{equation}\label{equ611}
\lim_{N\to\infty}\widetilde{E}_{i_N} \bigl(f\bigl(\widetilde
{Y}^N_{[N^2t]}\bigr) \bigr)
= \widetilde{E}_x(f(\widetilde{X}_t)) = \widetilde
{Q}_tf(x)
\end{equation}
for any $f\in C[0,1]$ with $f(1)=0$. If $i_N = i_N(x) = \min\{
[Nx]+1,1\}$,
the convergence is uniform in $x$, and is also uniform in $t$ for $0\le
t\le T$ for any $T>0$.
\end{theorem}

A weaker version of (\ref{equ611}) could be obtained from Theorem \ref{theo42} directly:
By (\ref{equ61})
%
\begin{eqnarray}\label{equ612}\qquad
\widetilde{E}_{i_N} (f(\widetilde{Y}^N_k) )
&=& \sum_{j=1}^{N-1} q_N^k(i_N,j)f\biggl(\frac jN\biggr)\nonumber \\
&=&\frac1{h_N(i_N)}
\sum_{j=1}^{N-1} p_N^k(i_N,j)h_N(j)f\biggl(\frac jN\biggr) \\
&=&\frac{s(x)}{h_N(i_N)} \frac1{s(x)}
\sum_{j=1}^{N-1} p_N^k(i_N,j)
\biggl[\frac{h_N(j)}{s(j/N)} s \biggl(\frac jN\biggr)
f \biggl(\frac jN\biggr) \biggr], \nonumber
\end{eqnarray}
where $1\le i_N\le N$ and $i_N/N\to x>0$. Now $g(x)=s(x)f(x)\in B_{01}$
if $f\in B_1$, and $\lim_{N\to\infty} h_N(i)/s(i/N)=s(1)$ uniformly for
$1\le i\le N$ by Lemma \ref{lemm53}. Hence by Theorem \ref{theo42}
\begin{eqnarray*}
\lim_{N\to\infty}\widetilde{E}_{i_N} \bigl(f\bigl(\widetilde
{Y}^N_{[N^2t]}\bigr) \bigr)
&=&\frac1{s(x)}Q_t(sf)(x) \\
&=&\frac1{s(x)}\int_0^1 p(t,x,y)s(y)f(y)m(dy)\\
&=&\widetilde{Q}_tf(x)
\end{eqnarray*}
uniformly for $s(x)\ge a$ for any $a>0$. However, this argument does not
extend to uniform convergence for $0\le x\le1$ nor to $x=0$.

\section{The limiting distribution of polymorphic sites}\label{sec7}
The purpose of this section is to prove Theorem \ref{theo21} and Corollary \ref{corr21} in
Section \ref{sec2}. By (\ref{equ23}), the expected value of the left-hand side of (\ref{equ213}) is
%
\begin{equation}\label{equ71}
E\Biggl(\sum_{i=1}^{N-1} f\biggl(\frac iN\biggr) N_{k_N}(i)\Biggr)
=
E\Biggl(\sum_{a=1}^{M_0} f\biggl(\frac{X^N_{1,a,k_N}}N\biggr)
+ \sum_{r=1}^{k_N}\sum_{b=1}^{M_r}
f\biggl(\frac{X^N_{2,b,r,k_N}}N\biggr) \Biggr).\hspace*{-42pt}
\end{equation}

\subsection[The first term in (7.1) (legacy polymorphisms)]{The first term in (\protect\ref{equ71}) (legacy polymorphisms)}\label{sec71}

The expected value of the sum in (\ref{equ71}) that corresponds to legacy
polymorphisms is
\begin{eqnarray*}
E\Biggl(\sum_{a=1}^{M_0} f\biggl(\frac{X^N_{1,a,k_N}}N\biggr)\Biggr)
&=&\sum_{i=1}^{N-1} E(N_0(i))
\sum_{j=1}^{N-1} p_N^{k_N}(i,j)f \biggl(\frac jN \biggr) \\
&=&\sum_{i=1}^{N-1} \omega^N_i
\sum_{j=1}^{N-1} p_N^{k_N}(i,j)f \biggl(\frac jN \biggr)\\
&=&\sum_{i=1}^{N-1} \omega^N_i Q^N_{k_N}f(i).
\end{eqnarray*}
By (\ref{equ61})
%
\begin{eqnarray}\label{equ72}\quad
\widetilde{Q}^N_kf(i)
&=&\sum_{j=1}^{N-1} q_N^k(i,j)f\biggl(\frac jN\biggr)
\nonumber\\
&=&\frac1{h_N(i)}
\sum_{j=1}^{N-1} p_N^k(i,j)h_N(j)f\biggl(\frac jN\biggr)
\\
&=&\frac1{h_N(i)} Q^N_k(h_Nf)(i), \nonumber
\end{eqnarray}
where $(h_Nf)(x)=h_N([Nx])f(x)$. Thus
%
\begin{eqnarray}\label{equ73}
\sum_{i=1}^{N-1} \omega^N_i Q^N_k f(i)
&=&\sum_{i=1}^{N-1} h_N(i)
\widetilde{Q}^N_k \biggl(\frac{f}{h_N} \biggr)(i) \omega^N_i
\nonumber\\[-8pt]\\[-8pt]
&=&
\sum_{i=1}^{N-1}
\biggl(\frac{h_N(i)}{s(i/N)} \biggr)
\widetilde{Q}^N_k \biggl(\frac f{h_N} \biggr)(i)
s \biggl(\frac{i}N \biggr) \omega^N_i. \nonumber
\end{eqnarray}
By Lemma \ref{lemm53}, $\lim_{N\to\infty} h_N(i_N)/s(i_N/N) = 1/s(1)$ uniformly
in $x$ for $i_N = i_N(x) = [Nx]$, and we can write
%
\begin{equation}\label{equ74}
\frac{f(i/N)}{h_N(i)} =g(i/N)
-g(i/N)\biggl(1 - \frac{s(i/N)/s(1)}{h_N(i)}\biggr)
\end{equation}
for $g(y)=s(1)f(y)/s(y)$. By the assumptions of Theorem \ref{theo21}, $g(y)$
extends to a continuous function on $C[0,1]$ with $g(1)=0$, so by
Theorem \ref{theo62}
%
\begin{equation}\label{equ75}
\lim_{N\to\infty} \widetilde{Q}^N_{[N^2t]} g([Nx]) = \widetilde
{Q}_tg(x)
\end{equation}
uniformly for $0\le x\le1$ and $0\le t\le T$ for any $T>0$. By (\ref{equ212})
%
\begin{equation}\label{equ76}
\lim_{N\to\infty} \sum_{i=1}^{N-1}
g \biggl(\frac iN \biggr) s \biggl(\frac iN \biggr) \omega^N_i
=\int_0^1 g(y)s(y)\nu(dy).
\end{equation}
Thus by (\ref{equ73}) and (\ref{equ75})
%
\begin{eqnarray}\label{equ77}
\lim_{N\to\infty} \sum_{i=0}^{N-1} \omega^N_i Q^N_k f(i)
&=& \int_0^1 \widetilde{Q}_t\bigl(s(1)f/s\bigr)(y)\frac
{s(y)}{s(1)}\nu(dy)\nonumber\\[-8pt]\\[-8pt]
&=& \int_0^1 Q_tf(y)\nu(dy).\nonumber
\end{eqnarray}
This completes the proof of Theorem \ref{theo21} for the legacy terms in (\ref{equ71}).

\subsection[The second term in (7.1) (new mutations)]{The second term in (\protect\ref{equ71}) (new mutations)}\label{sec72}
The expected value of the double sum in (\ref{equ71}), which corresponds to new
mutations, is
%
\begin{eqnarray}\label{equ78}\qquad
E\Biggl(\sum_{r=1}^{k_N}\sum_{b=1}^{M_r}
f\biggl(\frac{X^N_{2,b,r,k_N}}N\biggr) \Biggr)
&=&
\sum_{r=1}^{k_N} E(M_r)
E_1\biggl(f\biggl(\frac{X^N_{k_N-r}}N\biggr) \biggr)\nonumber\\
&=&
\mu_N \sum_{r=1}^{k_N}
Q^N_{k_N-r} f(1) \\
&=& \mu_N
\sum_{r=0}^{k_N-1} Q^N_r f(1), \nonumber
\end{eqnarray}
where $Q^N_r f(1)=\sum_{j=1}^{N-1} p_N(1,j)f(j/N)$. As in (\ref{equ72})
%
\begin{eqnarray}\label{equ79}\qquad
\mu_N\sum_{r=0}^{k_N-1} Q^N_r f(1)
&=&\mu_N
h_N(1)\sum_{r=0}^{k_N-1}
\widetilde{Q}^N_r(f/h_N)(1)\nonumber\\[-8pt]\\[-8pt]
&=&
(N\mu_N) (Nh_N(1))
\int_0^{k_N/N^2} \widetilde{Q}^N_{[N^2u]}(f/h_N)(1) \,du,
\nonumber
\end{eqnarray}
where $(f/h_N)(y)=f(y)/h_N([Ny]+1)$.

As $N\to\infty$, $N\mu_N\to\theta$ and $k_N/N^2\to t<\infty$ by (\ref{equ29})
and, since $s'(0)=1$, $Nh_N(1)\to1/s(1)$ by Lemma \ref{lemm53}. Thus by (\ref{equ79})
and (\ref{equ74}) and (\ref{equ75})
%
\begin{equation}\label{equ710}
\lim_{N\to\infty}
\mu_N \sum_{r=0}^{k_N-1} Q^N_r f(1)
=
\frac{\theta}{s(1)}\int_0^t \widetilde{Q}_ug(0) \,du,
\end{equation}
where $g(y)=s(1)f(y)/s(y)$. By (\ref{equ66}) and (\ref{equ69}), the Green operator
\begin{eqnarray*}
\widetilde{G}f(x) &=&\int_0^\infty\widetilde{Q}_tf(x)\,dt
=\int_0^1 \widetilde{g}(x,y)f(y)\widetilde{m}(dy) \\
&=&
\int_0^1 \biggl(\frac1{s(x\vee y)} - \frac1{s(1)}\biggr) f(y)\widetilde
{m}(dy)\nonumber
\end{eqnarray*}
is bounded in the supremum norm, and
\begin{eqnarray*}
&&
\frac{\theta}{s(1)}\int_0^t \widetilde{Q}_u g(0)\,du
\\
&&\qquad=\frac{\theta}{s(1)}
\biggl(\int_0^\infty\widetilde{Q}_u g(0)\,du - \int_t^\infty\widetilde
{Q}_u g(0)\,du\biggr)\\
&&\qquad=
\frac{\theta}{s(1)}\int_0^\infty\widetilde{Q}_u(g-\widetilde
{Q}_t g)(0)\,du\\
&&\qquad=
\frac{\theta}{s(1)}\int_0^1 \biggl(\frac1{s(y)}-\frac1{s(1)}\biggr)
\bigl(g(y) - \widetilde{Q}_t g(y)\bigr)\widetilde{m}(dy)
\\
&&\qquad=
\frac{\theta}{s(1)}\int_0^1 \frac{s(1)-s(y)}{s(1)s(y)}
\bigl(g(y) - \widetilde{Q}_t g(y)\bigr)s(y)^2 m(dy).
\end{eqnarray*}
Since $\widetilde{Q}_t g(y) = (1/s(y)) Q_t(sg)$ by (\ref{equ610})
and $g(y) =
s(1)f(y)/s(y)$, we conclude $g(y)-\widetilde{Q}_t g(y) =
(s(1)/s(y))(f(y)-Q_t f(y))$ and hence
\[
\frac{\theta}{s(1)}\int_0^t \widetilde{Q}_ug(0)\,du
=\theta\int_0^1 \frac{s(1)-s(y)}{s(1)-s(0)}
\bigl(f(y) - Q_tf(y)\bigr)m(dy).
\]
This is the second term on the right-hand side of (\ref{equ213}) and completes
the proof of Theorem \ref{theo21}.

\subsection[Proof of Corollary 2.1]{Proof of Corollary \protect\ref{corr21}}\label{sec73}

The first term $\int_0^1 Q_tf(y)\nu(dy)$ in the last line in
(\ref{equ213}) equals
%
\begin{equation}\label{equ711}
\int_0^1 \widetilde{Q}_t(f/s)(y) s(y)\nu(dy)
\end{equation}
as in (\ref{equ77}) or (\ref{equ610}), where $f(y)/s(y)$ is bounded and $\int_0^1
s(y)\nu(dy)<\infty$. Thus the integral in (\ref{equ711}) is $O(e^{-\lambda_1t})$
by (\ref{equ69}) and converges to zero as $t\to\infty$. By the same argument
applied with $(s(1)-s(x))m(dx)$ in place of $\nu(dx)$,
\[
\int_0^1 \frac{s(1)-s(x)}{s(1)-s(0)} Q_tf(x) m(dx)
=O(e^{-\lambda_1t})
\]
as $t\to\infty$ as well. This completes the proof of Corollary
\ref{corr21}.\vadjust{\goodbreak}

\section{The limiting numbers of fixations}\label{sec8}

\subsection[Proof of Theorem 2.2]{Proof of Theorem \protect\ref{theo22}}\label{sec81}
By Bartlett's theorem \cite{rKingman}, the number of process\-es
$\{X_{1,a,k}\}$ and $\{X_{2,b,r,k}\}$ that have been trapped at state
$N$ by time $k$ is Poisson with mean
%
\begin{equation}\label{equ81}
E (N_k(N) ) =
\sum_{i=1}^{N-1} \omega^N_i p_N^k(i,N)
+\mu_N \sum_{r=0}^{k-1} p_N^r(1,N).
\end{equation}
The first term on the right-hand side of (\ref{equ81}) corresponds to legacy
polymorphisms and the second to new polymorphisms. The proof of
Theorem \ref{theo22} for both terms depends on:
\begin{lem}\label{lemm81}
Let $i_N$ be integers such that $1\le i_N\le N$ and $i_N/N\to x$ for
some $x$ with $0\le x\le1$. Then
%
\begin{equation}\label{equ82}
\lim_{N\to\infty}P_{i_N}\biggl(\frac1{N^2} T^N_N \le t
\Bigm| T^N_N < T^N_0\biggr)
=P_x(T_1\le t\mid T_1<T_0).
\end{equation}
\end{lem}

We first show how Lemma \ref{lemm81} implies Theorem \ref{theo22}.
\begin{pf*}{Proof of Theorem \ref{theo22} given Lemma \ref{lemm81}}
The first sum in (\ref{equ81}) can be written
%
\begin{equation}\label{equ83}
\sum_{i=1}^{N-1} \omega^N_i p_N^k(i,N)
=\sum_{i=1}^{N-1}
\frac{p_N^k(i,N)}{s(i/N)} s \biggl(\frac iN\biggr) \omega^N_i.
\end{equation}
Let $i=i_N$ be integers with $1\le i_N\le N$ and $i/N\to x$, and assume
$k/N^2\to t$ as in Section \ref{sec2}. Then
\begin{eqnarray*}
\frac{p_N^k(i,N)}{s(i/N)}
&=& \frac{h_N(i)}{s(i/N)}
P_i(X^N_k=N \mid T^N_N<T^N_0) \\
&=& \biggl(\frac{h_N(i)}{s(i/N)}\biggr)
P_i \biggl(\frac1{N^2}T^N_N \le\frac{k}{N^2}
\Bigm| T^N_N<T^N_0 \biggr) \\
&\to& \frac1{s(1)}
P_x (T_1 \le t \mid T_1<T_0 )
\end{eqnarray*}
by Lemma \ref{lemm53} and (\ref{equ82}), with uniform convergence for $0\le x\le1$ if $i_N
= i_N(x) = \min\{[Nx]+1,1\}$. Thus by (\ref{equ83}) and (\ref{equ212})
%
\begin{eqnarray}\label{equ84}\qquad
\lim_{N\to\infty} \sum_{i=1}^{N-1} \omega^N_i p_N^k(i,N)
&=&\frac1{s(1)} \int_0^1 P_x (T_1\le t
\mid T_1<T_0 ) s(x)\nu(dx) \nonumber\\
&=&
\int_0^1 P_x(T_1\le t
\mid T_1<T_0) P_x(T_1<T_0) \nu(dx) \\
&=&
\int_0^1 P_x (T_1\le t) \nu(dx). \nonumber
\end{eqnarray}
The second term on the right-hand side of (\ref{equ81}) is
%
\begin{eqnarray}\label{equ85}\qquad
&&
\mu_N \sum_{r=0}^{k-1} p_N^r(1,N)
\nonumber\\
&&\qquad= \mu_N h_N(1) \sum_{r=0}^{k-1}
P_1(T^N_N \le r\mid T^N_N < T^N_0)
\nonumber\\[-8pt]\\[-8pt]
&&\qquad=
(N\mu_N)(Nh_N(1)) \int_0^{k/N^2}
P_1\biggl(\frac1{N^2} T^N_N \le\frac{[N^2u]}{N^2}
\Bigm| T^N_N < T^N_0\biggr) \,du \nonumber\\
&&\qquad\to
\frac\theta{s(1)} \int_0^t P_0(T_1\le u
\mid T_1<T_0 )\,du \nonumber
\end{eqnarray}
by Lemma \ref{lemm81}, since $N\mu_N\to\theta$ and $Nh_N(1)\to1/s(1)$ as in (\ref{equ79}).
Combining (\ref{equ84}) and (\ref{equ85}) completes the proof of Theorem \ref{theo22}.
\end{pf*}
\begin{pf*}{Proof of Lemma \ref{lemm81}}
Assume $1\le i_N\le N$ and $i_N/N\to x$ for $0\le x\le1$ as before.
Then by Theorem \ref{theo62}
\[
\lim_{N\to\infty} \widetilde{E}_{i_N}\bigl(f\bigl(Y^N_{[N^2t]}\bigr)\bigr)
=\widetilde{E}_x(f(X_t)),
\]
if $f\in C[0,1]$, $0\le f(x)\le1$, and $f(1)=0$. Since $N$ is a trap,
this implies
%
\begin{equation}\label{equ86}
\widetilde{P}_x(T_1 > t) \le \liminf_{N\to\infty}
\widetilde{P}_{i_N} \biggl(\frac1{N^2}T^N_N > t \biggr)
\end{equation}
and
\[
\widetilde{E}_x(T_1) \le \liminf_{N\to\infty}
\widetilde{E}_{i_N} \biggl(\frac1{N^2}T^N_N \biggr)
\]
by two applications
of Fatou's theorem. If we are able to prove
%
\begin{equation}\label{equ87}
\widetilde{E}_x(T_1) = \lim_{N\to\infty}
\widetilde{E}_{i_N} \biggl(\frac1{N^2}T^N_N \biggr) <
\infty,
\end{equation}
then a standard compactness argument for weak convergence would imply
equality in (\ref{equ86}) with $\liminf$ replaced by $\lim$, which would prove
Lemma \ref{lemm81}. Hence it is sufficient to prove (\ref{equ87}).

Consider an arbitrary birth-and-death Markov chain $X_n$ on the state
space $S=\{0,1,\ldots,N\}$ with absorbing endpoints. As in the Moran
model in (\ref{equ21}), assume that the transition function can be written
%
\begin{equation}\label{equ88}
p_{ij} =
\cases{
q_i\dvtx j=i+1, \cr
r_i\dvtx j=i, \cr
p_i\dvtx j=i-1,}
\end{equation}
where $p_{00}=p_{NN}=1$, $p_i+r_i+q_i=1$, and $p_i,q_i>0$ for $1\le i<N$
(see, e.g., \cite{rKarlin}, pages 50, 92--94). Then we have
the following lemma.
\begin{lem}\label{lemm82}
Let $g_{ij}=\sum_{n=0}^\infty p{}^{(n)}_{ij}$ where $p^{(n)}$
are the matrix powers of $p$ in (\ref{equ88}). Then
%
\begin{eqnarray}\label{equ89}
h_i &=& P_i(T_N<T_0) = \Biggl(\sum_{j=1}^{i-1}\alpha_j \Biggr)
\bigg/ \Biggl(\sum_{j=1}^{N-1}\alpha_j \Biggr)
=A_i/A_N, \nonumber\\[-8pt]\\[-8pt]
g_{ij} &=& A_N \frac{h_{i\wedge j}(1 - h_{i\vee j})}{\alpha_jq_j}
\nonumber
\end{eqnarray}
for $1\le i,j\le N-1$, $i\wedge j=\min\{i,j\}$ and $i\vee
j=\max\{i,j\}$. Here
%
\begin{equation}\label{equ810}
\alpha_j = \prod_{k=1}^j (p_k/q_k) \quad\mbox{and}\quad
A_i = \sum_{j=0}^{i-1} \alpha_j,\qquad (1\le i\le N).
\end{equation}
\end{lem}
\begin{pf}
The probabilities $h_i=P_i({T_N<T_0})$ satisfy the recurrence
$h_i = E_i(h_{X_1}) = p_i h_{i-1} + r_i h_i + q_i h_{i+1}$
for $0<i<N$. This implies $h_{i+1}-h_i=(p_i/q_i)(h_i-h_{i-1})$ and hence
$h_{i+1}-h_i = \alpha_i(h_1-h_0)$. By summation, $h_i-h_0 = A_i(h_1-h_0)$.
The conditions $h_0=0$ and $h_N=1$ imply $h_1=1/A_N$ and hence
$h_i=A_i/A_N$.

The Green matrix $g_{ij}$ satisfies the recurrence
\[
\sum_{j=0}^N p_{ir}g_{r,j}
=p_i g_{i-1,j} + r_i g_{ij} + q_i g_{i+1,j}
=g_{ij}-\delta_{ij}
\]
and thus $ g_{i+1,j}-g_{ij} = (p_i/q_i)(g_{i,j}-g_{i-1,j}) -\delta_{ij}/q_j$.
The boundary conditions $g_{0j}=g_{Nj}=0$ for $0<j<N$ and arguments
similar to those for $h_i$ lead to the formula for $g_{ij}$.
\end{pf}
%

For the Moran model (\ref{equ21}), $p_k/q_k=1/(1+\sigma_N)$,
$\alpha_j=(1+\sigma_N)^{-j}$, and $ A_i=\frac{1+\sigma_N}{\sigma_N}
(1 - (1+\sigma_N)^{-i})$. Thus by (\ref{equ89})
%
\begin{equation}\label{equ811}
h_N(i) =P_i(T^N_N<T^N_0) =\frac{A_i}{A_N} =
\frac{ 1 - (1+\sigma_N)^{-i}}{1 - (1+\sigma_N)^{-N}},
\end{equation}
which gives a second derivation of Lemma \ref{lemm51}.

If $f(i)=1$ for $0<i<N$ and $f(0)=f(N)=0$, then
%
\begin{eqnarray}\label{equ812}
E_i(T^N_N) &=& E_i\Biggl(\sum_{k=0}^{T^N_N-1} f(X_k)\Biggr)
= \sum_{k=0}^\infty E_i(f(X_k))\nonumber\\[-8pt]\\[-8pt]
&=& \sum_{j=1}^{N-1} g_N(i,j),\nonumber
\end{eqnarray}
where $g_N(i,j)=\sum_{k=0}^\infty p_N^k(i,j)$. Then by (\ref{equ812}) for the dual
Markov chain in (\ref{equ61}) and Lemma \ref{lemm82}
%
\begin{eqnarray}\label{equ813}\qquad
\widetilde{E}_i\biggl(\frac1{N^2}T^N_N\biggr)
&=&\frac1{N^2} \sum_{j=1}^{N-1} \widetilde{g}_N(i,j)
=\frac1{Nh_N(i)}\frac1{N}\sum_{j=1}^{N-1} g_N(i,j)h_N(j)
\nonumber\\
&=&\frac{A_N}{Nh_N(i)}\frac1{N}
\sum_{j=1}^{N-1} \bigl(h_N(i\wedge j)\bigl(1-h_N(i\vee j)
\bigr) \bigr)
\frac{h_N(j)}{\alpha_jq_j} \nonumber\\
&=&\frac{1 - (1+\sigma_N)^{-N}}{\sigma_N Nh_N(i)}
\\
&&{}\times\frac1N \sum_{j=1}^{N-1}
\biggl( \biggl(1 + \sigma_N\frac jN \biggr)
\bigl(h_N(i\wedge j)\bigl(1-h_N(i\vee j)\bigr) \bigr)\nonumber\\
&&\hspace*{120.2pt}{} \times
h_N(j) \frac{(1+\sigma_N)^j}{j/N(1-j/N)} \biggr)
\nonumber
\end{eqnarray}
by (\ref{equ21}) and (\ref{equ89}), since $\alpha_j=(1+\sigma_N)^{-j}$. By Lemma \ref{lemm53}
\[
\lim_{N\to\infty}
\frac{h_N(i)(1-h_N(j))}{s(i/N)(s(1)-s(j/N))}
=\frac1{s(1)^2}
\]
uniformly for $1\le i,j\le N-1$. Since $\sigma_N\sim\gamma/N$
by (\ref{equ29}), the
terms in the sum in (\ref{equ813}) are uniformly bounded. Then by Lemma \ref{lemm53}
again, with $i/N\to x$ and $j/N\to y$ in (\ref{equ813}),
%
\begin{eqnarray}\label{equ814}\qquad
&&\lim_{N\to\infty} \widetilde{E}_i\biggl(\frac1{N^2}T^N_N\biggr) \nonumber\\
&&\qquad =
\biggl(\frac{1 - e^{-\gamma}}{\gamma}\biggr)\frac1{s(x)}
\int_0^1 \frac{s(x\wedge y)(1-s(x\vee y))}{s(1)^2}
s(y) \frac{e^{\gamma y}}{y(1-y)}\,dy \\
&&\qquad =
\frac1{s(x)}
\int_0^1 g(x,y)s(y)m(dy) =
\int_0^1 \widetilde{g}(x,y) \widetilde{m}(dy) =\widetilde
{E}_x(T_1) \nonumber
\end{eqnarray}
for $g(x,y)$ in (\ref{equ43}), $m(dx)$ in (\ref{equ42}), $\widetilde{g}(x,y)$ in (\ref{equ68}) and
$\widetilde{m}(dx)$ in (\ref{equ65}). This completes the proof of (\ref{equ87}) and
hence of
Lemma \ref{lemm81} and Theorem \ref{theo22}.
\end{pf*}

\subsection[Proof of (2.17) after Theorem 2.2]{Proof of (\protect\ref{equ217}) after Theorem \protect\ref{theo22}}\label{sec82}
For fixed $t>0$
%
\begin{eqnarray}\label{equ815}
\widetilde{P}_x(T_1\le t) &=& 1 - \widetilde{P}_x(T_1>t) =1
- \int_0^1 q(t,x,y) \widetilde{m}(dy) \nonumber\\[-8pt]\\[-8pt]
&=& 1 - \int_0^1 \frac{p(t,x,y)}{s(x)} s(y) m(dy) \nonumber
\end{eqnarray}
by (\ref{equ68}). Thus
\begin{eqnarray*}
P_x(T_1\le t)
&=& P_x(T_1<T_0) P_x(T_1\le t\mid T_1<T_0) \\
&=& \bigl(s(x)/s(1)\bigr) \widetilde{P}_x(T_1\le t) \\
&=& \frac1{s(1)} \biggl( s(x) -\int_0^1 p(t,x,y) s(y) m(dy) \biggr).
\end{eqnarray*}
By (\ref{equ43}), (\ref{equ48}) and (\ref{equ46})
\[
\biggl|\frac{d\alpha_n(x)}{ds(x)} \biggr|
\le\lambda_n\int_0^1 |\alpha_n(y)|m(dy)
\le C_3\lambda_n^3
\]
and
%
\begin{equation}\label{equ816}
\frac{\partial}{\partial s(x)} p(t,x,y)
=\sum_{n=1}^{\infty}e^{-\lambda_nt} \frac{d\alpha
_n(x)}{ds(x)}\alpha_n(y)
\end{equation}
converges uniformly for $t\ge a>0$ and $0\le x,y\le1$. Thus by (\ref{equ815})
%
\begin{eqnarray}\label{equ817}
\widetilde{P}_0(T_1\le t)
&=& 1 - \lim_{a\to0} \int_0^1 \frac{p(t,a,y)}{s(a)} s(y) m(dy)
\nonumber\\
&=& 1 - \int_0^1 \frac{\partial}{\partial s(x)}
p(t,0+,y) s(y) m(dy) \\
&=& 1 - \int_0^1 q(t,0,y) s(y)^2 m(dy). \nonumber
\end{eqnarray}

%

%
\printaddresses

\end{document}